\documentclass[12pt]{article}
\usepackage{amsmath}
\usepackage{amssymb}
\usepackage{latexsym}
\usepackage{amssymb}
\usepackage{epsf}

\newtheorem{theorem}{Theorem}
\newtheorem{lemma}{Lemma}
\newtheorem{proposition}{Proposition}

\numberwithin{equation}{section}

\begin{document}
\date{}
\author{M.I.Belishev \thanks {Saint-Petersburg Department of
                 the Steklov Mathematical Institute (POMI), 27 Fontanka,
                 St. Petersburg 191023, Russia; belishev@pdmi.ras.ru.
                 Supported by RFBR grants 11-01-00407A and NSh-4210.2010.1.}}

\title{A unitary invariant of semi-bounded operator
in reconstruction of manifolds}
\maketitle
\begin{abstract}

With a densely defined symmetric semi-bounded operator of nonzero
defect indexes $L_0$ in a separable Hilbert space ${\cal H}$ we
associate a topological space $\Omega_{L_0}$ ({\it wave spectrum})
constructed from the reachable sets of a dynamical system governed
by the equation $u_{tt}+(L_0)^*u=0$. Wave spectra of unitary
equivalent operators are homeomorphic.

In inverse problems, one needs to recover a Riemannian manifold
$\Omega$ via dynamical or spectral boundary data. We show that for
a generic class of manifolds, $\Omega$ is {\it isometric} to the
wave spectrum $\Omega_{L_0}$ of the minimal Laplacian
$L_0=-\Delta|_{C^\infty_0(\Omega\backslash
\partial \Omega)}$ acting in ${\cal H}=L_2(\Omega)$, whereas $L_0$
is determined by the inverse data up to unitary equivalence.
Hence, the manifold can be recovered (up to isometry) by the
scheme `data $\Rightarrow L_0 \Rightarrow \Omega_{L_0}
\overset{\rm isom}= \Omega$'.

The wave spectrum is relevant to a wide class of dynamical
systems, which describe the finite speed wave propagation
processes. The paper elucidates the operator background of the
boundary control method (Belishev`1986), which is an approach to
inverse problems based on their relations to control theory.
\end{abstract}
\setcounter{section}{-1}

\section{Introduction}
\subsection{Motivation}
The paper introduces the notion of a {\it wave spectrum} of a
symmetric semi-bounded operator in a Hilbert space. The impact
comes from inverse problems of mathematical physics; the following
is one of the motivating questions.

Let $\Omega$ be a smooth compact Riemannian manifold with the
boundary $\Gamma$, $-\Delta$ the (scalar) Laplace operator,
$L_0=-\Delta|_{C^\infty_0(\Omega \backslash \Gamma)}$ the {\it
minimal Laplacian} in ${\cal H}=L_2(\Omega)$. Assume that we are
given with a unitary copy $\widetilde L_0=U L_0 U^*$ in a space
$\widetilde {\cal H} =U{\cal H}$ (but $\Omega, {\cal H}$ and $U$
are unknown!). To what extent does $\widetilde L_0$ determine the
manifold $\Omega$?

So, we have no points, boundaries, tensors, etc, whereas the only
thing given is an operator $\widetilde L_0$ in a Hilbert space
$\widetilde {\cal H}$. Provided the operator is unitarily
equivalent to $L_0$, is it possible to `extract' $\Omega$ from
$\widetilde L_0$? Such a question is an invariant version of
various setups of dynamical and spectral inverse problems on
manifolds \cite{BIP97}, \cite{BIP07}.

\subsection{Content}
Substantially, the answer is affirmative: for a generic class of
manifolds, any unitary copy of the minimal Laplacian determines
$\Omega$ up to isometry (Theorem 1). A wave spectrum is a
construction that realizes the determination $\widetilde L_0
\Rightarrow \Omega$ and thus solves inverse problems. In more
detail,
\medskip

\noindent $\bullet$\,\,With a closed densely defined symmetric
semi-bounded operator $L_0$ of nonzero defect indexes in a
separable Hilbert space ${\cal H}$ we associate a topological
space $\Omega_{L_0}$ (its wave spectrum). The space consists of
the atoms of a {\it lattice with inflation} determined by $L_0$.
The lattice is composed of reachable sets of an abstract {\it
dynamical system with boundary control} governed by the
evolutionary equation $u_{tt}+L_0^*u=0$. The wave spectrum is
endowed with a relevant topology.

Since the definition of $\Omega_{L_0}$ is of invariant character,
the spectra $\Omega_{L_0}$ and $\Omega_{\widetilde L_0}$ of
unitarily equivalent operators $L_0$ and $\widetilde L_0$ turn out
to be homeomorphic. So, a wave spectrum is a (hopefully, new)
unitary invariant of a symmetric semi-bounded operator.
\medskip

\noindent $\bullet$\,\, A wide generic class of the so-called {\it
simple manifolds} is introduced\footnote{Roughly speaking, a
simplicity means that the symmetry group of $\Omega$ is trivial.}.
The central Theorem 1 establishes that for a simple $\Omega$, the
wave spectrum of its minimal Laplacian $L_0$ is metrizable and
isometric to $\Omega$. Hence, any unitary copy $\widetilde L_0$ of
$L_0$ determines the simple $\Omega$ up to isometry by the scheme
$\widetilde L_0 \Rightarrow \Omega_{\widetilde L_0}\overset{\rm
isom}=\Omega_{L_0} \overset{\rm isom}= \Omega$. In applications,
it is the procedure, which recovers manifolds by the {\it boundary
control method} \cite{BIP97},\cite{BIP07}: concrete inverse data
determine the relevant $\widetilde L_0$, what enables one to
realize the scheme.
\medskip

\noindent $\bullet$\,\, We discuss one more option: elements of
the space ${\cal H}$  can be realized as `functions' on
$\Omega_{L_0}$ \footnote{In the BC-method, such an option is
interpreted as {\it visualization of waves} \cite{BIP07}.}.
Hopefully, this observation can be driven at a functional model of
a class of $L_0$s and/or Spaces of Boundary Values. Presumably,
this model will be {\it local}, i.e., satisfying ${\rm supp\,}
(L_0^{\rm mod})^*y \subseteq {\rm supp\,} y$.

\subsection{Comments}
$\bullet$\,\,The concept of wave spectrum summarizes rich
`experimental material' accumulated in inverse problems of
mathematical physics in the framework of the BC-method, and
elucidates its operator background. For the first time,
$\Omega_{L_0}$ has appeared in \cite{BKac89} in connection with
the M.Kac problem; its later version (called a {\it wave model})
is presented in \cite{BIP07} (sec 2.3.4). Owing to its invariant
nature, $\Omega_{L_0}$ promises to be useful for further
applications to unsolved inverse problems of elasticity theory,
electrodynamics, graphs, etc.

Our paper is of pronounced interdisciplinary character. `Wave'
terminology, which we use, is motivated by close relations to
applications.
\medskip

\noindent $\bullet$\,\, The path from $L_0$ to $\Omega_{L_0}$
passes through an intermediate object, which is a sublattice of
the lattice ${\mathfrak L}({\cal H})$ of subspaces of the space ${\cal H}$. Section 1
is an excursus to the lattice theory, in course of which we
introduce {\it lattices with inflation}. The wave spectrum appears
as a set of {\it atoms} of the relevant lattice with inflation
determined by $L_0$.
\medskip

\noindent $\bullet$\,\,We give attention to connections of our
approach with C*-algebras. As is shown, if $\Omega$ is a compact
manifold then $\Omega_{L_0}$ is identical to the Gelfand spectrum
of the algebra of continuous functions $C(\Omega)$. By the recent
trend in the BC-method, to recover unknown manifolds via boundary
inverse data is to find spectra of relevant algebras determined by
the data \cite{BSobolev}. We hope for utility and further
promotion of this trend.
\medskip

\noindent{\bf Acknowledgements}\,\,\, I would like to thank
A.B.Alexandrov for kind consultations, and V.I.Vasyunin and N.Wada
for useful discussions on the subject of the paper.

\section{Lattices with inflation}
Reducing the volume of the paper, we do not prove Propositions.
The proofs are quite elementary and typical technique is
demonstrated in Appendix.

\subsection{Basic objects}
{\bf 1.\,Lattice.}\,\,\,Let $\mathfrak L$ be a {\it lattice}, i.e.
a partially ordered set ({\it poset}) with the order $\leqslant$
and operations $a \wedge b=\inf \{a, b\}$, $a \vee b =\sup
\{a,b\}$. Also, we assume that ${\mathfrak L}$ is endowed with the
{\it least element} $0$ satisfying $0<a$ for $a \not= 0$
\cite{Birk}.

The {\it order topology} on ${\mathfrak L}$ is introduced through
the {\it order convergence}: $x_j \to x$ if there are the nets
$\{a_j\}_{j \in J}\! \uparrow$ and $\{b_j\}_{j \in J}\!\downarrow$
($J$ is a directed set) such that $a_j \leqslant x_j \leqslant
b_j$ and $\sup \{a_j\}=x=\inf \{b_j\}$ holds (we write
$a_j\uparrow x$ and $b_j\downarrow x$). For an $A \subset
{\mathfrak L}$, the inclusion $x \in \overline A$
occurs\footnote{Everywhere $\overline{(\,\,)}$ denotes a
topological closure. In some places, to avoid the confusion, we
specify the space.} if and only if there are $a_j, b_j \in A$ such
that $a_j\uparrow x$ and/or $b_j\downarrow x$ \cite{Birk}.
\smallskip

{\bf Example 1.}\,\,\,The lattice ${\mathfrak L}=2^\Omega$ of
subsets of a set $\Omega$ with the order $\leqslant = \subseteq$,
operations $\wedge=\cap,\,\,\vee=\cup$, and $0=\emptyset$.

{\bf Example 2.}\,\,\,The (sub)lattice ${\cal O} \subset 2^\Omega$ of
open sets of a topological space $\Omega$. The convergence
$\omega_j\uparrow \omega$ means $\omega=\sup \{\omega_j\}=\cup_j\,
\omega_j$. The convergence $\omega_j\downarrow \omega$ means
$\omega=\inf \{\omega_j\}={\rm int} \cap_j\, \omega_j$, where ${\rm int}
A$ is the set of interior points of $A \subset \Omega$.
\bigskip

\noindent{\bf 2.\,Inflation.}\,\,\,For a lattice ${\mathfrak L}$, the set
${\cal F}_{\mathfrak L} := {{\cal F}} \left([0,\infty); {\mathfrak L} \right)$ of ${\mathfrak L}$-valued
functions is also a topologized lattice with respect to the
point-wise order, operations, and topology.

{\bf Definition 1.}\,\,\, A map $I: {\mathfrak L} \to {\cal F}_{\mathfrak L}$ is said to be
an {\it inflation} if for all $a, b \in {\mathfrak L}$ and $s,t \in
[0,\infty)$ one has

\noindent(i)\,\,\, $(Ia)(0)=a$ and $I 0_{\mathfrak L}\,=\, 0_{{\cal F}_{\mathfrak L}}$

\noindent(ii)\,\,\, $a \leqslant b$ and $s \leqslant t$ imply
$(Ia)(s) \leqslant (Ib)(t)$.

Inflation is injective: $I^{-1}f=f(0)$ on $I{\mathfrak L}$.
\smallskip

{\bf Example 3.}\,\,\,$\Omega$ is a metric space with the distance
${\rm d}$. For a subset $A \subset \Omega$, by $A^t:=\{x \in
\Omega\}\,|\,\,{\rm d}(x,A)<t\}$ ($t>0$) we denote its metric
neighborhood, ant put $A^0:=A,\,\,\,\emptyset^t=\emptyset$. The
map $M: 2^\Omega \to {\cal F}_{2^\Omega},\,\,(MA)(t):=
A^t,\,\,t\geqslant 0$ is a {\it metric inflation}. The image
$M2^\Omega$ is a semilattice: $Ma \vee Mb = M(a \vee b) \in
M2^\Omega$. The image of open sets  is a (sub)semilattice $M{\cal
O} \subset {\cal F}_{{\cal O}} \subset {\cal F}_{2^\Omega}$.
\bigskip

\noindent{\bf 3. Atoms. Basic construction.}\,\,\,Let ${\cal P}$
be a poset with the least element $0$. An $\alpha \in {\cal P}$ is
called an {\it atom} if $0<a\leqslant \alpha$ implies $a=\alpha$
\cite{Birk}. By ${\rm At\,} {\cal P}$ we denote the set of atoms.
\smallskip

{\bf Example 4.}\,\,\,Each atom of $2^\Omega$ is a single point
set: ${\rm At\,} 2^\Omega=\left\{\{x\}\,|\,\, x \in \Omega\right\}$.

{\bf Example 5.}\,\,\,If the open sets of a topological space
$\Omega$ are infinitely divisible\footnote{i.e., for any
$\emptyset \not=A \in {\cal O}$ there is $\emptyset \not=B \in {\cal O}$
such that $B \subset A$ and $A\backslash B\not=\emptyset$.} then
${\rm At\,} {\cal O}=\emptyset$.
\smallskip

Inflation preserves atoms: $I {\rm At\,} {\mathfrak L} \subseteq {\rm At\,} I{\mathfrak L}$.

For any lattice with inflation, the set $\Omega_{I{\cal L}}:={\rm At\,}
\overline{I{\mathfrak L}} \subset {\cal F}_{\mathfrak L}$ (the closure in topology of ${\cal F}_{\mathfrak L}$)
is well defined \footnote{But the case $\Omega_{I{\cal L}}=\emptyset$
is not excluded.}. This set is a key object of the paper. Namely,
the following effect will be exploited: there are lattices and
inflations such that ${\rm At\,} {\mathfrak L}=\emptyset$ but ${\rm At\,}
\overline{I{\mathfrak L}}\not=\emptyset$. {\bf Inflation can create atoms!}
\smallskip

There is a natural topology on $\Omega_{I{\cal L}}\subset {\cal
F}_{\mathfrak L}$. For atoms $\alpha, \beta \in \Omega_{I{\cal
L}}$, we say that $\alpha$ influences on $\beta$ at the moment $t$
if $\alpha(t)\wedge \beta(\varepsilon)\not=0_{\mathfrak L}$ for
any positive $\varepsilon$. Define $t_{\alpha \beta}:=\inf
\{t\geqslant 0\,|\,\,\alpha(t)\wedge
\beta(\varepsilon)\not=0_{\mathfrak L} \,\,\forall
\varepsilon>0\}$. If $\alpha(t)\wedge
\beta(\varepsilon)=0_{\mathfrak L}$ for all positive $t$ and
$\varepsilon$, we put $t_{\alpha \beta}=\infty$.

A function $\tau_{I{\mathfrak L}}: \Omega_{I{\cal L}} \times \Omega_{I{\cal L}} \to
[0,\infty) \cup\{\infty\}, \,\tau_{I{\mathfrak L}}(\alpha,
\beta):=\max\{t_{\alpha \beta}, t_{\beta \alpha}\}$ is called an
{\it interaction time}.

Define the `balls' $B^r[\alpha]:=\{\beta \in
\Omega_{I{\cal L}}\,|\,\,\tau_{I{\mathfrak L}}(\alpha, \beta)<r\}\,
(r>0),\,\,B^0[\alpha]:=\alpha$.

{\bf Definition 2.}\,\, By $(\Omega_{I{\cal L}}, \tau_{I{\mathfrak L}})$ we
denote the topological space that is the set $\Omega_{I{\cal L}}$
endowed with the minimal topology, which contains all balls.

Surely, at this level of generality, to expect for rich properties
of this space is hardly reasonable. However, in `good' cases the
function $\tau_{I{\mathfrak L}}$ turns out to be a metric.
\begin{proposition}\label{P1} Let $(\Omega, {\rm d})$ be a complete metric
space, ${\mathfrak L}=2^\Omega$, $I=M$ (see Example 3). The
correspondence $\Omega \ni x \leftrightarrow M\{x\} \in {\cal
F}_{2^\Omega}$ is a bijection between the sets $\Omega$ and
$\Omega_{M2^\Omega}= {\rm At\,} \overline{M2^\Omega}={\rm At\,}
M2^\Omega= M{\rm At\,} 2^\Omega=\{M\{x\}\,|\,\,x \in \Omega\}$.
The equality $\tau_{M2^\Omega}(M\{x\},M\{y\})={\rm d} (x,y)$
holds. Function $\tau_{M2^\Omega}$ is a metric on atoms, whereas
$(\Omega_{M2^\Omega},\tau_{M2^\Omega})$ is a metric space
isometric to $(\Omega, {\rm d})$. The isometry is realized by the
bijection $M\{x\} \leftrightarrow x$.
\end{proposition}
\bigskip

There are another topologies on atoms, which are also inspired by
the metric topology. The first one is introduced via closure
operation: for a set $W \subset \Omega_{I\cal L}$, we put
$$\overline W:=\left\{\alpha \in \Omega_{I\cal
L}\,\,\biggl|\,\,\bigvee \limits_{\beta \in W} \beta \overset{\cal
F}\geqslant \alpha\right\}.$$ It is easy to check that the map $W
\mapsto \overline W$ satisfies all Kuratovsky's axioms and, hence,
determines a unique topology $\rho_{I\cal L}$ in $\Omega_{I\cal
L}$. Note a certain resemblance (duality) of such a topology to
Jacobson's topology on the set $\cal I$ of primitive ideals of a
C*-algebra $\cal A$. Namely, for a $W \subset \cal I$, one defines
its closure by $$\overline W:=\left\{i \in {\cal
I}\,\,\biggl|\,\,\bigcap \limits_{b \in W} b \subseteq i\right\}$$
(see, e.g., \cite{Mur}).

One more topology is the following. We define the `balls' by
$$B^r[\alpha]:=\left\{\beta \in \Omega_{I\cal L}\,\,\biggl|\,\,\exists t_0=t_0(\alpha, \beta, r)>0
\,\,\,{\rm s.t.}\,\,\,0 \not=\beta(t_0) \overset{\cal L}\leqslant
\alpha(r)\right\} \quad (r>0).$$ As one can verify, the system
$\{B^r[\alpha]\}_{\alpha \in \Omega_{I\cal L}, \,r>0}$ is a base
of topology. Hence, it determines a unique topology that we denote
by $\sigma_{I\cal L}$.
\medskip

If ${\cal L}=2^{{\mathbb R}^n}$ and $I$ is the (Euclidean) metric
inflation, the topologies $\tau_{I\cal L},\, \rho_{I\cal L}$, and
$\sigma_{I\cal L}$ coincide with the standard Euclidean metric
topology in ${\mathbb R}^n$.
\bigskip

{\bf 4. Isomorphic lattices.}\,\,Let ${\mathfrak L}$ and ${\mathfrak L}'$ be two
lattices with inflations $I$ and $I'$ respectively. We say them to
be {\it isomorphic} through a bijection $i:{\mathfrak L}\to{\mathfrak L}'$ if $i$
preserves the order, lattice operations, and $i(IA)=I'i(A)$ holds
for all $A \in {\mathfrak L}$.

The bijection $i$ is extended to bijection on functions $i:{\cal
F}_{\mathfrak L} \to {\cal F}_{{\mathfrak L}'}$ by
$(if)(t):=i(f(t)),\,\,t\geqslant 0$. The following fact is quite
obvious.
\begin{proposition}\label{P2} If the lattices with
inflation ${\mathfrak L}$ and ${\mathfrak L}'$ are isomorphic then
the spaces $(\Omega_{I{\cal L}}, \tau_{I{\mathfrak L}})$ and
$(\Omega_{I'{\cal L}'}, \tau_{I'{\mathfrak L}'})$ are
homeomorphic. The homeomorphism is realized by the bijection $i$
on atoms.
\end{proposition}

\subsection{Lattices in metric space}
{\bf 5.\, Lattice ${\cal O}$.}\,\,\,Return to Example 3 and assume in
addition that

\noindent{\bf A1.}\,\,$\Omega$ is a complete metric space

\noindent{\bf A2.}\,\,All the balls $\{x\}^t$ are compact and
$\{x\}^t \backslash \{x\}^s\not=\emptyset$ as $s<t$.

By A2, open sets are infinitely divisible. Therefore, we have
${\rm At\,} {\cal O}=\emptyset$.
\smallskip

Fix an $x \in \Omega$ and define the functions $x_*, x^* \in {\cal
F}_{{\cal O}}$: $x_*(t):=\{x\}^t$ as $t>0$, $x_*(0):=0_{\cal O},$
and $x^*(t):={\rm int} \overline{\{x\}^t}\,\,$ as $t\geqslant 0$.
Evidently, we have $x_*\leqslant x^*$ in ${\cal F}_{\cal O}$. The
upper function satisfies $x^*=\lim_{\varepsilon \to
0}\,M\left(\{x\}^\varepsilon\right) \in \overline{M{\cal
O}},\,\,x^*(0)=0_{\cal O}$. The `clearance' between the functions
is small: $\overline{x_*(t)}=\overline{x^*(t)},\,\,t\geqslant 0$.

Since $x^* \in \overline{M{\cal O}}$, the segment $[x_*, x^*]:=\{f
\in {\cal F}_{\cal O}\,|\,\,x_* \leqslant f \leqslant x^*\}$
intersects with $\overline{M{\cal O}}$. The poset $[x_*, x^*]\cap
\overline{M{\cal O}}$ is a closed subset in ${\cal F}_{\cal O}$
bounded from below. Hence, it contains minimal elements, which can
be easily recognized as the atoms of $\overline{M{\cal O}}$. So,
$\Omega_{M{\cal O}}:={\rm At} \overline{M{\cal O}}\not=
\emptyset$.

{\bf Example 6.}\,\,\,For $\Omega \subseteq {\mathbb R}^n$ one has
$x_*=x^*$. Therefore, each segment $[x_*,x^*]$ contains one (and
only one) atom $\{x\}^t, \, t\geqslant 0$. We don't know whether
the same is correct for a Riemannian manifold $\Omega$.
\smallskip

For an atom $\alpha \in {\rm At} \overline{M{\cal O}}$, define a {\it kernel}
$\dot \alpha:= \cap_{t>0}\,\alpha(t) \subset \Omega$.
\begin{proposition}\label{P3} For each $\alpha$, its kernel $\dot
\alpha$ consists of a single point $x_\alpha \in \Omega$. Each
atom $\alpha$ belongs to the segment $[(x_\alpha)_*,
(x_\alpha)^*]$. If $\dot \alpha=\dot \beta$ then
$\overline{\alpha(t)}=\overline{\beta(t)},\,\, t\geqslant 0$
holds.
\end{proposition}
These facts follow from a general lemma stated and proved in
Appendix.
\smallskip

With each $x \in \Omega$ one associates the class of atoms
$\langle \alpha \rangle_x:=[x_*, x^*]\cap {\rm At}
\overline{M{\cal O}}$. For $\alpha, \beta \in \langle \alpha
\rangle_x$ one has $\overline{\alpha(t)}=\overline{\beta(t)}\,
(=\overline{\{x\}^t}\,),\,\,t\geqslant 0$. Hence, $\alpha$ and
$\beta$ interact at any $t>0$. As a result, we have $\tau_{M{\cal
O}}(\alpha, \beta)=0$.

The relation $\{\alpha \sim \beta\}\Leftrightarrow
\{\tau_{M{\cal O}}(\alpha, \beta)=0\}$ is an equivalence on
$\Omega_{M{\cal O}}$. The factor-set
${\Omega^\prime_{M{\cal O}}}:=\Omega_{M{\cal O}}/{\sim}$ is bijective to
$\Omega$ through the map $\langle \alpha \rangle \mapsto \dot \alpha$. The
function $\tau^\prime_{M{\cal O}}(\langle \alpha \rangle, \langle \beta \rangle):=
\tau_{M{\cal O}}(\alpha, \beta)$ is a metric on $\Omega^\prime_{M{\cal O}}$.
The equality $\tau^\prime_{M{\cal O}}(\langle \alpha \rangle_x, \langle \alpha
\rangle_y)={\rm d}(x,y)$ is valid for all $x,y \in \Omega$ and we
conclude the following.
\begin{proposition}\label{P4}
The metric space $(\Omega^\prime_{M{\cal O}},\tau^\prime_{M{\cal
O}})$ is isometric to $(\Omega, {\rm d})$. The isometry is
realized by the bijection $\langle \alpha \rangle_x
\leftrightarrow x$.
\end{proposition}
\bigskip

{\bf 6.\,Lattice ${{\cal O}^{\rm reg}}$.}\,\,For a set $A \subset \Omega$, denote
by $\partial A:=\overline A \cap \overline{\Omega \backslash A}$ its boundary. Note
that $\partial (A \cap B) \subseteq \partial A \cup \partial B$ and $\partial (A \cup
B) \subseteq \partial A \cup \partial B$. It is convenient to put $\partial
\Omega=\partial \emptyset=\emptyset$.

Recall that we deal with complete and locally compact metric
spaces (see A1,2). In addition, assume that $\Omega$ is endowed
with a Borel measure $\mu$ such that

\noindent{\bf A3.}\,\,For {\it any} $A \subset \Omega$ and $t>0$,
the relation $\mu(\partial A^t)=0$ holds.

{\bf Example 7.}\,\,$\Omega$ is a smooth Riemannian manifold with
the canonical measure (volume). In particular, $\Omega \subseteq
{\mathbb R}^n$ with the Lebesgue measure \cite{Fed}.
\smallskip

{\bf Definition 3.}\,\,An open set $A \subset \Omega$ is called
{\it regular} \footnote{Our definition is similar to (but differs
from) the definition of regularity in \cite{Birk}, p 216.} if
$\mu(\partial A)=0$. The system of regular sets is denoted by ${{\cal O}^{\rm reg}}$.

As is evident, ${{\cal O}^{\rm reg}}$ is a sublattice in ${\cal
O}$. It is a base of ${\cal O}$: each open set is a sum of regular
sets (balls). By A3, ${{\cal O}^{\rm reg}}$ is invariant w.r.t.
the metric inflation: $(M {{\cal O}^{\rm reg}})(t)\subset {{\cal
O}^{\rm reg}},\, t\geqslant 0$. In other words, we have $M {{\cal
O}^{\rm reg}} \subset {\cal F}_{{{\cal O}^{\rm reg}}}$.
\smallskip

Fix an $x \in \Omega$. Note that $x_*, x^* \in {\cal F}_{{\cal O}^{\rm reg}}$ and $x^*
\in \overline{M{{\cal O}^{\rm reg}}}$. Using the arguments quite analogous to ones, which
have led to Proposition \ref{P4}, and factorizing the set of atoms
w.r.t. the same equivalence $\sim$, one can arrive at the
following result.
\begin{proposition}\label{P5}
The metric space $(\Omega^\prime_{M{{\cal O}^{\rm
reg}}},\tau^\prime_{M{{\cal O}^{\rm reg}}})$ is isometric to
$(\Omega, {\rm d})$. The isometry is realized by the bijection
$\langle \alpha \rangle_x \leftrightarrow x$.
\end{proposition}

The operation $A \mapsto A^*:={\rm int} (\Omega \backslash A)$ is well
defined on ${{\cal O}^{\rm reg}}$ and called a {\it pseudo-complement} \cite{Birk}.
The relations $A \cap A^*=\emptyset$ and $A \subseteq (A^*)^*$ are
valid.
\bigskip

{\bf 7.\,Lattice ${\mathfrak R}$.}\,\,\,Introduce an equivalence on ${{\cal O}^{\rm reg}}$:
we put $A \simeq B$ if $\overline A = \overline B$. Define ${\mathfrak R}:={{\cal O}^{\rm reg}}
/\simeq\,$. By $[A]$ we denote the equivalence class of $A$.

Endow ${\mathfrak R}$ with the order and operations:

\noindent $[A]\leqslant[B]$ if $A \subseteq$ B

\noindent $[A]\wedge [B]:=[A \cap B]$, $[A]\vee [B]:=[A \cup B]$

\noindent $[A]^\bot:=[A^*]\,\,\,\,\left(=[{\rm int} (\Omega \backslash
A)]\right)$.

\noindent The least and greatest elements are $0:=[\emptyset]$ and
$1:=[\Omega]$.

One can easily check the well-posedness of these definitions and
prove the following relations:

\noindent$[A]\wedge [A]^\bot=0,\,\,[A]\vee [A]^\bot=1$

\noindent$\left([A]\wedge [B]\right)^\bot=[A]^\bot\vee
[B]^\bot,\,\,\left([A]\vee [B]\right)^\bot=[A]^\bot\wedge [B]^\bot$.

\noindent Hence ${\mathfrak R}$ is a lattice with the {\it complement}
$[\,\cdot\,]^\bot$ \,\cite{Birk}.
\smallskip

For $f \in {\cal F}_{{{\cal O}^{\rm reg}}}$, define $[f]\in {\cal
F}_{\cal R},\,\,[f](t):=[f(t)],\,\,t\geqslant 0$.

Introduce the {\it metric inflation} on ${\mathfrak R}$ by $M:
{\mathfrak R} \to {\cal F}_{\mathfrak
R},\,\,\left(M[A]\right)(t):=\left[(MA\right)(t)]=[A^t], \,\,t
\geqslant 0$.

The relation ${\rm At} \overline{M{\cal
R}}=\{[\alpha]\,|\,\,\alpha \in {\rm At}\overline{M{{\cal O}^{\rm
reg}}}\}$ holds. The map $A\mapsto[A]$ identifies the atoms
belonging to the same class: if $\alpha, \beta \in \langle\alpha
\rangle_x$ then
$\overline{\alpha(t)}=\overline{\beta(t)},\,\,t\geqslant 0$ that
implies $[\alpha]=[\beta]$. By this, the set $\Omega_{M{\mathfrak
R}}={\rm At} \overline{M{\mathfrak R}}$ is bijective to $\Omega$,
whereas the `interaction time' $\tau_{M{\mathfrak R}}$ turns out
to be a metric.
\begin{proposition}\label{P6}
The metric space $(\Omega_{M{\mathfrak R}},\tau_{M{\mathfrak R}})$
is isometric to $(\Omega, {\rm d})$. The isometry is realized by
the bijection $[\alpha]\leftrightarrow x_\alpha$.
\end{proposition}
\bigskip

{\bf 8.\,Lattice ${{\mathfrak R}^{\cal H}}$.}\,\,\,Introduce a
Hilbert space ${\cal H}:=L_{2, \mu}(\Omega)$.

For a measurable set $A \subset \Omega$, define the subspace ${\cal H}
A:=\{\chi_A y\,|\,\, y \in {\cal H}\}$, where $\chi_A$ is the indicator
of $A$. Such subspaces are called {\it geometric}. If $A\in {{\cal O}^{\rm reg}}$
then $\mu(\overline A \backslash A)= \mu (\partial A)=0$ that leads to ${\cal H}
\overline A={\cal H} A$.

{\bf Definition 4.}\,\,If $A \in {{\cal O}^{\rm reg}}$, we say the subspace ${\cal H} A$
to be {\it regular}. The system of regular subspaces is denoted by
${{\mathfrak R}^{\cal H}}$.
\smallskip

Let ${{\mathfrak L}(\cal H)}$ be the lattice of subspaces of the space ${\cal H}$ (see item
10 below). The system ${{\mathfrak R}^{\cal H}} \subset {{\mathfrak L}(\cal H)}$ is a sublattice.

Introduce a map $i: {\mathfrak R} \to {{\mathfrak R}^{\cal H}},\,\,i[A]:={\cal H} A$. As is easy to
check, it preserves the operations and complement\footnote{The
latter means $i([A]^\bot)=({\cal H} A)^\bot = {\cal H} \ominus {\cal H} A={\cal H}
A^*$.}.

Extend $i$ to functions: for an $f \in {\cal F}_{\mathfrak R}$ we
put $if \in {\cal F}_{{\mathfrak R}^{\cal H}} \subset {\cal
F}_{{\mathfrak L}(\cal H)},\,\, (if)(t):=i(f(t)), \,\,t \geqslant
0$. Also, define a {\it metric inflation} on ${{\mathfrak R}^{\cal
H}}$ by $iM: {{\mathfrak R}^{\cal H}} \to {\cal F}_{{\mathfrak
L}(\cal H)},\,\,\newline \left(iM {\cal H} A\right)(t):={\cal
H}\left(\left(MA\right)(t)\right) =i[A^t]={\cal H} A^t, \,\,t
\geqslant 0$.

 Thus, $i$ is an isomorphism of
lattices with inflation. Propositions \ref{P3} and \ref{P6} lead
to the following result.
\begin{proposition}\label{P7}
The metric space $(\Omega_{iM{{\mathfrak R}^{\cal
H}}},\tau_{iM{{\mathfrak R}^{\cal H}}})$ is isometric to $(\Omega,
{\rm d})$. The isometry is realized by the bijection
$i[\alpha]\leftrightarrow x_\alpha$.
\end{proposition}

The meaning of the passage ${{\cal O}^{\rm reg}} \to {{\mathfrak R}^{\cal H}}$ is that it `codes'
regular sets in Hilbert terms. Later in inverse problems, we will
determine the Hilbert lattices from inverse data, and then
`decode' them, i.e., extract information about geometry of
$\Omega$.
\bigskip

{\bf 9.\,Dense sublattice.}\,\,We say a system of subsets ${\cal
N} \subset {{\cal O}^{\rm reg}}$ to be {\it dense in} ${{\cal
O}^{\rm reg}}$, if for any $x \in \Omega$ and $A \in {{\cal
O}^{\rm reg}},\,\,x \in A$ there is an $N \in {\cal N}$ such that
$x \in N \subset A$. If, moreover, ${\cal N}$ is a sublattice such
that $M{\cal N}\subseteq {\cal N}$ holds, we call it a {\it dense
$M$-invariant sublattice in} ${{\cal O}^{\rm reg}}$.

Let ${\mathfrak R}_{\cal N} \subseteq {\mathfrak R}$ be the image of ${\cal N}$ through the map
$A \mapsto [A]$. The following fact can be derived as a
consequence of density.
\begin{proposition}\label{P8}
If ${\cal N}$ is a dense $M$-invariant sublattice, then the metric
space $(\Omega_{M {\mathfrak R}_{\cal N}},\tau_{M {\mathfrak
R}_{\cal N}})$ is isometric to $(\Omega, {\rm d})$. The isometry
is realized by the bijection $[\alpha]\leftrightarrow x_\alpha$.
\end{proposition}
\medskip

Let ${\mathfrak R}^{\cal H}_{\cal N} \subset {\mathfrak R}^{\cal H} \subset {{\mathfrak L}(\cal H)}$ be the image of
${\cal N}$ through the map ${\cal N} \ni N \mapsto {\cal H} N \in {\mathfrak R}^{\cal H}$. The
image is an $iM$-invariant sublattice in ${\mathfrak R}^{\cal H}$. The next
result is a straightforward consequence of the previous one.
\begin{proposition}\label{P9}
If ${\cal N} \subset {{\cal O}^{\rm reg}}$ is a dense
$M$-invariant sublattice then the metric space
$(\Omega_{iM{\mathfrak R}^{\cal H}_{\cal N}},\tau_{iM{\mathfrak
R}^{\cal H}_{\cal N}})$ is isometric to $(\Omega, {\rm d})$. The
isometry is realized by the bijection $i[\alpha]\leftrightarrow
x_\alpha$.
\end{proposition}
Later, in applications, we will deal with concrete $\Omega$ and
${\cal N}$.
\medskip

The relation $(\Omega, {\rm d})\overset{\rm
isom}=(\Omega_{iM{\mathfrak R}^{\cal H}_{\cal N}},\tau_{iM{\mathfrak R}^{\cal H}_{\cal N}})$ is the final
goal of our excursus to the lattice theory. It represents the
original metric space as collection of atoms of relevant Hilbert
lattice with inflation. This representation will play the key role
in reconstruction of $\Omega$ via inverse data.

\section{Wave spectrum}
\subsection{Inflation in ${\cal H}$}

{\bf 10.\,Basic objects.}\, Let ${\cal H}$ be a separable Hilbert
space, ${{\mathfrak L}(\cal H)}$ the lattice of its (closed)
subspaces equipped with the order $\leqslant = \subseteq$,
operations ${\cal A} \wedge {\cal B}={\cal A} \cap {\cal B}$,
${\cal A} \vee {\cal B}=\overline{\{a+b\,|\,\,a \in {\cal A}, b
\in {\cal B}\}}$, the complement ${\cal A} \mapsto {\cal
A}^\bot={\cal H} \ominus {\cal A}$, and extremal elements
$0=\{0\}$, $1={\cal H}$. A {\it sublattice} of ${{\mathfrak
L}(\cal H)}$ is its subset closed w.r.t. the operations and
complement. Each sublattice contains $0$ and $1$.

By $P_{\cal A}$ we denote the (orthogonal) projection onto ${\cal
A} \in {{\mathfrak L}(\cal H)}$. Also, if ${\cal A}$ is a non-closed lineal set, we
put $P_{\cal A}:=P_{\overline{\cal A}}$.

Let $\mathfrak B ({\cal H})$  be the algebra of bounded operators.
For an ${\cal S} \subseteq {{\mathfrak B}(\cal H)}$, by ${\rm Proj} {\cal S}$ we denote the set of
projections belonging to ${\cal S}$.

For a lattice ${\mathfrak L} \subseteq {{\mathfrak L}(\cal H)}$,
with a slight abuse of notation, we put ${\rm Proj} {\mathfrak L}
:=\{P_{\cal A}\,|\,\,{\cal A} \in {\mathfrak L}\}$. The map
${\mathfrak L} \ni {\cal A} \mapsto P_{\cal A} \in {\rm Proj}
{\mathfrak L}$ induces the lattice structure on ${\rm Proj}
{\mathfrak L}$: $P_{\cal A} \wedge P_{\cal B} = P_{{\cal A} \wedge
{\cal B}},\,\,P_{\cal A} \vee P_{\cal B} = P_{{\cal A} \vee {\cal
B}},\,\,(P_{\cal A})^\bot=P_{{\cal A}^\bot}$. For $P,Q \in {\rm
Proj} {\mathfrak L}$ the relation $P\leqslant Q$ means ${\rm Ran}
P \subseteq {\rm Ran} Q$ and holds if and only if
$(Px,x)\leqslant(Qx,x), \,\, x \in {\cal H}$. The extremal
elements of ${\rm Proj} {\mathfrak L}$ are the zero and unit
operators $\mathbb O$ and $\mathbb I$.

The same map relates the order topology on ${{\mathfrak L}(\cal H)}$ with the strong
operator topology on ${{\mathfrak B}(\cal H)}$: ${\cal A}=\lim {\cal A}_j$ in ${{\mathfrak L}(\cal H)}$ if and only
if $P_{\cal A}=s\!-\!\lim P_{{\cal A}_j}$ in ${{\mathfrak B}(\cal H)}$.
\medskip

{\bf Convention.}\,\,The metric inflation $iM$ in ${\cal H}=L_{2,
\mu}(\Omega)$ is defined on a sublattice ${\mathfrak R}^{\cal H} \subset {{\mathfrak L}(\cal H)}$
(item 8). In contrast to it, in what follows we deal with
inflations {\it defined on the whole} ${{\mathfrak L}(\cal H)}$.

For an inflation $I:{{\mathfrak L}(\cal H)} \to {\cal
F}_{{\mathfrak L}(\cal H)}$, we denote ${\cal A}^t:= (I{\cal
A})(t),\,\,t\geqslant 0$. It is also convenient to regard
inflation as an operation on projections: $I:{\rm Proj}
{{\mathfrak B}(\cal H)} \to {\cal F}_{{\rm Proj} {{\mathfrak
B}(\cal H)}},\,(IP)(t)=P^t:=P_{(I {\rm Ran} P)(t)},\,\,t \geqslant
0$.
\smallskip

A lattice ${\mathfrak L} \subset {{\mathfrak L}(\cal H)}$ is said
to be {\it $I$-invariant} if $I{\mathfrak L} \subset {\cal
F}_{\mathfrak L}$ holds, i.e. ${\cal A} \in {\mathfrak L}$ implies
$(I{\cal A})(t)\in {\mathfrak L},\,\,t \geqslant 0$.

{\bf Definition 5.}\,\, Let ${\mathfrak f} \subseteq {{\mathfrak L}(\cal H)}$ be a family of
subspaces. Define ${\mathfrak L}[I,{\mathfrak f}]\subseteq {{\mathfrak L}(\cal H)}$ as the minimal
$I$-invariant lattice, which contains ${\mathfrak f}$.
\bigskip

\noindent{\bf 11.\,Spectra.}\, Let ${\cal H}$ and $I$ be given, ${\mathfrak L}$
be an $I$-invariant lattice. Recall that the space of atoms with
the `interaction' topology was introduced in item 3.

{\bf Definition 6.}\,\,The space $\Omega^{\rm At}_{\mathfrak L}:=(\Omega_{I {\mathfrak L}},
\tau_{I {\mathfrak L}})$ is called an {\it atomic spectrum} of the lattice
${\mathfrak L}$.
\smallskip

There is a version of this notion. Each function $f \in
\overline{I{\mathfrak L}} \subset {\cal F}_{\mathfrak L}$ is an
increasing family of subspaces $\{f(t)\}_{t\geqslant 0} \subset
{{\mathfrak L}(\cal H)},\,\,$, i.e. a {\it nest} \cite{Don}. The
corresponding nest of projections $\{P^t_f\}_{t \geqslant
0},\,\,P^t_f:=P_{f(t)}$ determines a self-adjoint operator
$E_f:=\int_0^\infty t\,dP^t_f$. It acts in ${\cal H}$ and is
called an {\it eikonal}. The set of eikonals is ${\rm Eik}
{\mathfrak L}:=\{E_f\,|\,\,f \in \overline{I{\mathfrak L}}\}$.

{\bf Definition 7.}\,\,A metric space
$\Omega^{\rm nest}_{\mathfrak L}:=\{E_\alpha\,|\,\,\alpha \in \Omega_ {I{\mathfrak L}}\}$
with the distance $\|E_\alpha - E_\beta\|$ is called a {\it nest
spectrum} of the lattice ${\mathfrak L}$.

{\bf Caution!}\,\,\,We do not assume $E_\alpha$ to be a bounded
operators, so that the `pathologic' situation ${\rm
dist}(E_\alpha, E_\beta)\equiv \infty$ is not excluded. However, a
`good' case, when all the differences $E_\alpha - E_\beta$ are
{\it bounded} operators, is realized in applications.
\medskip

One more version is the following.

Let us say that we deal with a {\it bounded case} if the set of
eikonals of the lattice is uniformly bounded: $\sup\{\|E\|\,|\,\,E
\in {\rm Eik}{\mathfrak L}\}<\infty$.

With a lattice ${\mathfrak L}$ one associates the von Neumann operator
algebra\footnote{i.e., a unital weekly closed self-adjoint
subalgebra of $\mathfrak B ({\cal H})$: see \cite{Mur}.} ${\mathfrak N}_{\mathfrak L}
\subseteq {{\mathfrak B}(\cal H)}$ generated by the projections of ${\mathfrak L}$, i.e., the
minimal von Neumann algebra satisfying ${\rm Proj} {\mathfrak L} \subseteq {\rm Proj}
{\mathfrak N}_{\mathfrak L}$.

In the bounded case, we have ${\overline {{\rm Eik} {\mathfrak L}}}^s\subset {\mathfrak N}_{\mathfrak L}$
(the closure in the strong operator topology). The elements of
this closure are also called {\it eikonals}.

The set ${\overline {{\rm Eik} {\mathfrak L}}}^s$ is partially
ordered: for two eikonals $E,\,E'$, we write $E \leqslant E'$ if
$(Ex,x) \leqslant (E'x,x),\,\,x \in {\cal H}$. An eikonal $E$ is
{\it maximal} if $E \leqslant E'$ implies $E=E'$. By $\Omega^{\rm
eik}_{\mathfrak L} \subset {\overline {{\rm Eik} {\mathfrak
L}}}^s$ we denote the set of maximal eikonals.

\begin{lemma}\label{L1}
In the bounded case, the set $\Omega^{\rm eik}_{\mathfrak L}$ is
nonempty.
\end{lemma}
{\bf Proof.}\,\,\,By the boundedness, any totally ordered family
of eikonals $\{E_j\}$ has an upper bound $s$-$\overline{\lim}
E_j$, which is also an eikonal. Hence, the Zorn lemma implies
$\Omega^{\rm eik}_{\mathfrak L} \not= \emptyset$.\,\,\,$\square$
\smallskip

{\bf Definition 8.}\,\,A metric space $\Omega^{\rm eik}_{\mathfrak L}$ with the
distance $\|E-E'\|$ is called an {\it eikonal spectrum} of the
lattice ${\mathfrak L}$.
\smallskip

In the general (unbounded) case, one can regularize the eikonals
as $E^\epsilon_f:=\int_0^\infty\frac{t}{1+\epsilon
t}\,dP^t_f\,\,\,(\epsilon>0)$ and deal with the corresponding
spectra ${\Omega}^{{\rm eik}, \epsilon}_{\mathfrak
L}\not=\emptyset$.
\bigskip

{\bf Remark.}\,\,\,Returning to Definition 6, one more option is
to define the atomic spectrum as $(\Omega_{I {\mathfrak L}},
\rho_{I {\mathfrak L}})$ or $(\Omega_{I {\mathfrak L}}, \sigma_{I
{\mathfrak L}})$ (see sec 1.1, item 3). Our reserve of concrete
examples is rather poor and provides no preferable choice.

\subsection{Inflation $I_L$}
{\bf 12.\,Dynamical system.}\,Let $L$ be a semi-bounded
self-adjoint operator in ${\cal H}$. Without lack of generality, we
assume that it is positive definite:
$$L=L^*\,=\,\int^\infty_\varkappa \lambda\,dQ_\lambda\,; \qquad (Ly,y)\geqslant
\varkappa \,\|y\|^2,\,\,\,y \in {\rm Dom\,} L \subset {\cal H}\,,
$$ where $d Q_\lambda$ is the spectral measure of $L$, $\varkappa>0$
is a constant.

Operator $L$ governs the evolution of a dynamical system
\begin{align}
\label{dynsystem_1} & v_{tt} + L v\,=\,h\,, \qquad t>0\\
\label{dynsystem_2} & v|_{t=0}\,=\,v_t|_{t=0}\,=\,0\,,
\end{align}
where $h \in L_2^{\rm loc}\left((0,\infty); {\cal H}\right)$ is a
${\cal H}$-valued function of time ({\it control}). Its solution
$v=v^h(t)$ is represented by the Duhamel formula
\begin{align}
\notag &  v^h(t)\,=\,\int_0^t L^{-\frac{1}{2}} \sin\left[(t-s)
L^{\frac{1}{2}}\right]\,h(s)\,ds =\\
\label{trajectory}& =\,\int^t_0 ds \int^\infty_0 \frac{\sin {\sqrt
\lambda}(t-s)}{\sqrt \lambda}\,d Q_\lambda\, h(s)\,, \qquad t
\geqslant 0\,
\end{align}
(see, e.g., \cite{BSol})\footnote{For $\varkappa\leqslant 0$,
problem (\ref{dynsystem_1}), (\ref{dynsystem_2}) is also well
defined but the representation (\ref{trajectory}) is of slightly
more complicated form.}. In system theory, $v^h$ is referred to as
a {\it trajectory}; $v^h(t)\in {\cal H}$ is a {\it state} at the
moment $t$. In applications, $v^h$ describes a {\it wave}
initiated by a source $h$.
\smallskip

Fix a subspace ${\cal A} \subseteq {\cal H}$. The set ${\cal
V}_{{\cal A}}^t:=\{v^h(t)\,|\,\,h \in L_2^{\rm loc}\left((0,
\infty); {\cal A}\right)\}$ of all states produced by ${\cal
A}$-valued controls is called {\it reachable} (at the moment $t$,
from the subspace $\cal A$). Reachable sets increase: ${\cal A}
\subseteq {\cal B}$ and $s\leqslant t$ imply ${\cal V}_{{\cal
A}}^s \subseteq  {\cal V}_{{\cal B}}^t$.
\bigskip

\noindent{\bf 13.\,Dynamical inflation.}\, With the system
(\ref{dynsystem_1}), (\ref{dynsystem_2}) one associates a map
$I_L: {{\mathfrak L}(\cal H)} \to {\cal F}_{{\mathfrak L}(\cal H)},
\,\,(I_L{\cal A})(0):={\cal A},\,(I_L{\cal A})(t):=\overline{{\cal V}_{{\cal
A}}^t},\,\,t>0$.
\begin{lemma}\label{L2}
$I_L$ is an inflation.
\end{lemma}
{\bf Proof.}\,\,\,The relation $(I_L{\cal A})(s)\subseteq
(I_L{\cal B})(t)$ as ${\cal A} \subseteq {\cal B}$ and
$0<s\leqslant t$ is a consequence of the general properties of
reachable sets. The only fact we need to verify is that the map
extends subspaces: ${\cal A} \subseteq (I_L{\cal A})(t),\,t>0$.
\smallskip

By $\chi_{[a,b]}$ we denote the indicator of the segment
$[a,b]\subset \mathbb R$. Fix an $r>0$ and $\varepsilon \in (0,r)$.
Define the  functions
$\varphi_\varepsilon(t):=\varepsilon^{-2}\chi_{[-\varepsilon,\varepsilon]}(t){\rm sign}(-t)$
and $\varphi_\varepsilon^r(t):=\varphi_\varepsilon(t-r+\varepsilon)$ for $t \in \mathbb
R$. Note that $\int_0^r\varphi_\varepsilon^r(t)\,f(t)\,dt \to -
f'(r)$ as $\varepsilon \to 0$ for a smooth $f$, i.e.,
$\varphi_\varepsilon^r(t)$ converges to $\delta^\prime(t-r)$ as a
distribution.

For $\lambda>0$, define a function
$$\psi_\varepsilon(\lambda):=\int_0^r
\frac{\sin[\sqrt{\lambda}\,(r-t)]}{\sqrt
\lambda}\,\varphi_\varepsilon^r(t)\,dt\,=\,\frac{2
\cos(\sqrt{\lambda} \,\varepsilon)-\cos(\sqrt{\lambda} \,2
\varepsilon) - 1}{\varepsilon^2 \lambda}.$$ Note that
$\psi_\varepsilon(\lambda) \underset{\varepsilon \to 0} \to 1$ as
$\varepsilon \to 0$ uniformly w.r.t. $\lambda$ in any segment
$[\varkappa, N]$.

Take a nonzero $y \in {\cal A}$ and consider (\ref{dynsystem_1}),
(\ref{dynsystem_2}) with the control
$h_\varepsilon(t)=\varphi_\varepsilon^r(t)\,y$. By the properties
of $\psi_\varepsilon$ one has
\begin{align*}
&\left\|y-v^{h_\varepsilon}(r)\right\|^2 = \langle{\rm see\,}
(\ref{trajectory})\rangle=\left\|y-\int_0^r dt
\int_{\varkappa}^\infty \frac{\sin[\sqrt{\lambda}\,(r-t)]}{\sqrt
\lambda}\,dQ_\lambda
[\varphi_\varepsilon(t) y]\right\|^2 =\\
& \left\|y - \int_{\varkappa}^\infty
\psi_\varepsilon(\lambda)\,dQ_\lambda y\right\|^2
=\left\|\int_{\varkappa}^\infty
\left[1-\psi_\varepsilon(\lambda)\right]\,dQ_\lambda
y\right\|^2=\\
& \int_{\varkappa}^\infty |1-\psi_\varepsilon(\lambda)|^2\,d
\|Q_\lambda y\|^2\,\underset{\varepsilon \to 0}\to 0.
\end{align*}
 The order of integration
change is easily justified by the Fubini Theorem.
\smallskip

Thus, $y=\underset{\varepsilon \to 0}\lim\, v^{h_\varepsilon}(r)$,
whereas $v^{h_\varepsilon}(r)\in (I_L{\cal A})(r)$ holds. Since
$(I_L{\cal A})(r)$ is closed in ${\cal H}$, we get $y \in (I_L{\cal
A})(r)$. Hence, ${\cal A} \subseteq (I_L{\cal A})(r),\,r>0$. \,\,\,$\square$
\medskip

So, each positive definite operator $L$ determines the inflation
$I_L$, which we call a {\it dynamical inflation}.

\subsection{Space $\Omega_{L_0}$.}
{\bf 14.\,Lattice ${\cal L}_{L, {\cal D}}$ and spectra.}\,\,Fix a subspace
${\cal D} \in {{\mathfrak L}(\cal H)}$ and say it to be a {\it directional subspace}.

Return to the system (\ref{dynsystem_1})--(\ref{dynsystem_2}).
Introduce the class ${\cal M}_{{\cal D}}:=\newline \left\{h \in
C^\infty\left([0, \infty); {\cal D}\right)\,|\,\,{\rm supp\,} h \subset
(0, \infty) \right\}$ of smooth ${\cal D}$-valued controls vanishing
near $t=0$. This class determines the sets
\begin{align}
\notag  &{\cal U}^t_{\cal D}:=\left\{h(t)- v^{h^{\prime
\prime}}(t)\,\biggl |\,h \in {\cal M}_{{\cal D}}\right\}\,=
\langle\,{\rm see\,} (\ref{trajectory})\, \rangle\,=\\
\label{reachable_sets_U} & \biggl\{h(t)-\int_0^t L^{-\frac{1}{2}}
\sin\left[(t-s) L^{\frac{1}{2}}\right]\, h^{\prime
\prime}(s)\,ds\,\biggl |\,\,h \in {\cal M}_{{\cal D}}\biggr\},
\qquad t \geqslant 0,
\end{align}
where $(\,\cdot\,)^\prime :=\frac{d}{dt}$. These sets are also
called {\it reachable}. As one can show, the sets ${\cal U}^t_{\cal D}$
increase as $t$ grows.

{\bf Definition 9.}\,\,The family of subspaces ${\mathfrak u}_{L,
{\cal D}}=\{\overline{{\cal U}^t_{\cal D}}\}_{t \geqslant
0}\subseteq {{\mathfrak L}(\cal H)}$ is called a {\it boundary
nest}.
\smallskip

The boundary nest determines the lattice ${\mathfrak L}_{L,
{\cal D}}:={\mathfrak L}[I_L,{\mathfrak u}_{L, {\cal D}}]$, which is the minimal $I_L$-invariant
sublattice in ${{\mathfrak L}(\cal H)}$ containing ${\mathfrak u}_{L, {\cal D}}$ (item 10).

The lattice determines the spectra $\Omega^{\rm At}_{{\mathfrak L}_{L, {\cal D}}}$ and
$\Omega^{\rm nest}_{{\mathfrak L}_{L, {\cal D}}}$. In the bounded case, the spectrum
$\Omega^{\rm eik}_{{\mathfrak L}_{L, {\cal D}}}$ is also well defined (item 11).
\bigskip

{\bf 15.\,Lattice ${\mathfrak L}_{L_0}$ and spectra.}\,\,Let $L_0$
be a closed densely defined symmetric semi-bounded operator with
{\it nonzero} defect indexes $n_{\pm} = n \leqslant \infty$. As is
easy to see, such an operator is necessarily unbounded. For the
sake of simplicity, we assume it to be positive definite: $(L_0 y,
y)\geqslant \varkappa \|y\|^2,\,\,y \in {\rm Dom} L_0$ with
$\varkappa>0$.

Let $L$ be the {\it Friedrichs extension} of $L_0$, so that $L=L^*
\geqslant \varkappa\,\mathbb I$ and $L_0 \subset L \subset L_0^*$
holds \cite{BSol}. Also, note that $1 \leqslant {\rm dim\,}{\rm
Ker\,}L_0^*= n \leqslant \infty$.
\smallskip

With the operator $L_0$ one associates two objects: the inflation
$I_L$ and the directional subspace ${\cal D}={\rm Ker}{L^*_0}$.
This pair determines the boundary nest ${\mathfrak
u}_{L_0}:={\mathfrak u}_{L, {\rm Ker}{L^*_0}}=\{\overline{{\cal
U}^t_{{\rm Ker} L_0}}\}_{t \geqslant 0}$ and the lattice
${\mathfrak L}_{L_0}:={\mathfrak L}_{L, {\rm Ker}{L^*_0}}$. The
nest and lattice determine the corresponding spectra, and we
arrive at the key subject of the paper.
\smallskip

{\bf Definition 10.}\,\,The space
$\Omega_{L_0}:=\Omega^{\rm At}_{{\mathfrak L}_{L_0}}$ is called a {\it wave
spectrum} of the (symmetric semi-bounded) operator $L_0$.

Recall that $\Omega_{L_0}$ is endowed with the `interaction time'
topology (item 3).
\smallskip

By analogy with the latter definition, one can introduce the
metric spaces $\Omega^{\rm nest}_{L_0}:=\Omega^{\rm nest}_{{\mathfrak L}_{L_0}}$ and,
in the bounded case, $\Omega^{\rm eik}_{L_0}:=\Omega^{\rm eik}_{{\mathfrak L}_{L_0}}$,
which are also determined by $L_0$.
\smallskip

As is evident from their definitions, the spectra are unitary
invariants of the operator.
\begin{proposition}\label{P10} If $U: {\cal H} \to
\tilde {\cal H}$ is a unitary operator and ${\tilde L}_0= U L_0 U^*$
then $\Omega_{{\tilde L}_0}$ is homeomorphic to $\Omega_{L_0}$. If
${\cal H}={\cal H}^1 \oplus {\cal H}^2$ and $L_0=L_0^1 \oplus L_0^2$ then
$\Omega_{L_0}=\Omega_{L_0^1} \cup \Omega_{L_0^2}$.
\end{proposition}
These properties motivate the use of term `spectrum'. The same
properties occur for $\Omega^{\rm nest}_{{\mathfrak L}_{L_0}}$ and
$\Omega^{\rm eik}_{L_0}$, replacing `homeomorphic' with `isometric'.
\bigskip

{\bf 16.\,Structures on $\Omega_{L_0}$.}\,\, The boundary nest
${\mathfrak u}_{L_0}$ can be regarded as an element (function) of the space
${\cal F}_{{\mathfrak L}(\cal H)}$. As such, it can be compared with the atoms, which
constitute the wave spectrum $\Omega_{L_0} \subset {\cal F}_{{\mathfrak L}(\cal H)}$.

{\bf Definition 11.}\,\,The set $\partial \Omega_{L_0}:=\{\alpha
\in \Omega_{L_0}\,|\,\,\alpha \leqslant {\mathfrak u}_{L_0}
\}$\footnote{Recall that $\alpha \leqslant {\mathfrak u}_{L_0}$ in
${\cal F}_{{\mathfrak L}(\cal H)}$ means that $\alpha(t)\subseteq
\overline{{\cal U}^t}$ holds for $t\geqslant 0$.} is said to be
the {\it boundary} of $\Omega_{L_0}$.

Also, it is natural to put $\partial\Omega^{\rm
nest}_{L_0}:=\partial \Omega_{L_0}$. In the bounded case, one
introduces the {\it boundary eikonal} $E^\partial=\int_0^\infty
t\,dP_{{\cal U}^t_{{\rm Ker} L_0}}$ and defines $\partial
\Omega^{\rm eik}_{L_0}=\{E\in \Omega^{\rm
eik}_{L_0}\,|\,\,E\geqslant E^\partial\}$ (see \cite{BArXiv}).
\medskip

There is a way to represent elements of ${\cal H}$ as `functions' on
the wave spectrum.

Fix an atom $\alpha \in \Omega_{L_0}:\,
\alpha=\alpha(t),\,\,t\geqslant 0$. Let
$P^t_\alpha:=P_{\alpha(t)}$ be the corresponding projections. For
$f,g \in {\cal H}$, we put $f \overset{\alpha}= g$ if there is
$\varepsilon=\varepsilon(f,g,\alpha)>0$ such that $P^t_\alpha
f=P^t_\alpha g$ as $t<\varepsilon$. The relation
$\overset{\alpha}=$ is an equivalence. The equivalence class
$[f]_\alpha=:G^f(\alpha)$ is called a {\it wave germ} (of the
element $f$ at the atom $\alpha$).

{\bf Definition 12.}\,\,The germ-valued function $G^f: \alpha
\mapsto [f]_\alpha,\,\,\alpha \in \Omega_{L_0}$ is called a {\it
wave image} of the element $f$.

The collection ${\cal G}:=\{G^f\,|\,\,f \in {\cal H}\}$ is a linear space
w.r.t. the point-wise algebraic operations: $\left(\lambda G^f+\mu
G^g\right)(\alpha):=[\lambda f+\mu g]_\alpha,\,\,\alpha \in
\Omega_{L_0}$. The linear map ${\cal I}:{\cal H} \ni f \mapsto G^f \in {\cal G}$
is called an {\it image operator}.

\section{DSBC}
\subsection{Green system}
{\bf 17. Ryzhov's axioms.} Consider a collection $\{ {\cal
H},{\cal B}; A,\Gamma_0, \Gamma_1\}$ of separable Hilbert spaces
${\cal H}$ and ${\cal B}$, and densely defined operators $A:{\cal
H} \to {\cal H}$ and $\Gamma_k:{\cal H} \to {\cal B}\,\,\,(k=0,1)$
connected via the {\it Green formula}
$$(Au,v)_{\cal H} - (u,Av)_{\cal H} =
(\Gamma_0 u,\Gamma_1 v)_{\cal B} - (\Gamma_1 u, \Gamma_0 v)_{\cal
B}\,.$$ The space ${\cal H}$ is called an {\it inner space}; ${\cal B}$
and $\Gamma_k$ are referred to as a {\it boundary values space}
and the {\it boundary operators} respectively \cite{Koch}. Such a
collection is said to be a {\it Green system}.

The following additional conditions are imposed.

\noindent{\bf R1.}\,\,${\rm Dom\,} \Gamma_k \supseteq {\rm Dom\,}
A$ holds. The restriction $A|_{{\rm Ker\,} \Gamma_0 \cap {\rm
Ker\,} \Gamma_1}=:L_0$ is a densely defined symmetric po\-si\-tive
definite operator with nonzero defect indexes. The relation
$\overline A=L_0^*$ is valid ('bar' is the operator closure).
\smallskip

\noindent{\bf R2.}\,\,The restriction $A|_{{\rm Ker\,} \Gamma_0}
=: L$ coincides with the Friedrichs extension of $L_0$, so that we
have $L_0 \subset L \subset L_0^*=\overline A$. Operator $L^{-1}$
is bounded and defined on ${\cal H}$.
\smallskip

\noindent{\bf R3.}\,\,The subspaces ${\cal A}:={\rm Ker\,} A$ and
${\cal D}:={\rm Ker\,} L_0^*$ are such that the relations
$\overline{\cal A} = {\cal D}$ and $\overline{\Gamma_0 {\cal A}}={\cal B}$
hold.
\smallskip

These conditions were introduced by V.A.Ryzhov \cite{Ryzh}, which
puts them as basic axioms. Note, that there are a few versions of
such an axiomatics but the one proposed in \cite{Ryzh} is most
relevant for applications to forward and inverse {\it
multidimensional} problems of mathematical physics.

The following consequences are derived from R1-4 \cite{Ryzh}.

\noindent{\bf C1.}\,\,The operator $\Pi:=\left(\Gamma_1
L^{-1}\right)^*: {\cal B} \to {\cal H}$ is bounded. The set ${\rm
Ran\,} \Pi$ is dense in ${\cal D}$.

\noindent{\bf C2.}\,\,The representation ${\cal A}=\{y \in {\rm
Dom\,} A\,|\,\,\Pi \Gamma_0 y = y\}$ is valid.

\noindent{\bf C3.}\,\,Since $L$ is the extension of $L_0$ by
Friedrichs, the relations ${\rm Dom\,} L_0 = L^{-1}[{\cal H}
\ominus {\cal D}]$ and $L_0\,=\,L|_{L^{-1}[{\cal H} \ominus {\cal
D}]}$ easily follows from the definition of such an extension (see
\cite{BSol}).
\bigskip

{\bf 18. Illustration.}\,\,Let $\Omega$ be a $C^\infty$-smooth
compact Riemannian manifold with the boundary $\Gamma$, $\Delta$
the (scalar) Beltrami-Laplace operator in ${\cal H}:=L_2(\Omega)$,
$\nu$ the outward normal on $\Gamma$, ${\cal B}:=L_2(\Gamma)$.
Denote \footnote{$H^k$ are the Sobolev classes; $H^2_0(\Omega)=\{y
\in H^2(\Omega)\,|\,y=|\nabla y|=0\,\, {\rm on}\,\Gamma\}$;
$\partial_\nu$ is the differentiation w.r.t. the outward normal on
$\Gamma$.} $A=-\Delta|_{H^2(\Omega)},\,\,\,\Gamma_0 u = u|_\Gamma,
\,\,\,\Gamma_1 u=\partial_\nu u|_\Gamma\,$, so that $\Gamma_{0,
1}$ are the {\it trace operators}. The collection $\{ {\cal
H},{\cal B}; A,\Gamma_0, \Gamma_1\}$ is a Green system. Other
operators, which enter in Ryzhov's axiomatics, are the following:
\smallskip

\noindent $\bullet$\,\,\,$L_0=-\Delta|_{H^2_0(\Omega)}$ is the
minimal Laplacian that coincides with the closure of
$-\Delta|_{C^\infty_0(\Omega\backslash
\partial \Omega)}$
\smallskip

\noindent $\bullet$\,\,\,$L=-\Delta|_{H^2(\Omega)\cap
H^1_0(\Omega)}$ is the self-adjoint Dirichlet Laplacian
\smallskip

\noindent $\bullet$\,\,\,$L_0^*=-\Delta |_{\{y \in {\cal
H}\,|\,\Delta y \in {\cal H}\}}$ is the maximal Laplacian
\smallskip

\noindent $\bullet$\,\,\,${\cal A} = \{y \in
H^2(\Omega)\,|\,\,\Delta y =0\}$ is the set of harmonic functions
of the class $H^2(\Omega)$
\smallskip

\noindent $\bullet$\,\,\,${\cal D}=\{y \in {\cal H}\,|\,\,\Delta
y=0\}$ is the subspace of all harmonic functions in $L_2(\Omega)$
\smallskip

\noindent $\bullet$\,\,\,$\Pi: {\cal B} \to {\cal H}$ is the
harmonic continuation operator (the Dirichlet problem solver):
$\Pi \varphi=u$ is equivalent to $\Delta u=0\,\,\,{\rm
in}\,\,\,\Omega,\,\,\,u|_{\Gamma}=\varphi$.

\subsection{Evolutionary DSBC}
{\bf 19. Dynamical system.}\,\,The Green system determines an
evolutionary {\it dynamical system with boundary control}
\begin{align}
\label{0.1} & u_{tt}+A\,u=0 &&{\rm in}\,\, {\cal H}, \quad 0<t<\infty   \\
\label{0.2} & u|_{t=0}=u_t|_{t=0}=0 &&{\rm in}\,\, {\cal H} \\
\label{0.3} & \Gamma_0 u=f(t) &&{\rm in}\,\,{\cal B}, \quad
0\leqslant t < \infty,
\end{align}
where $f$ is a {\it boundary control}, $u=u^f(t)$ is the solution
({\it wave}). The space of controls ${\cal F}=L^{\rm
loc}_2\left((0,\infty); \cal B\right)$ is said to be {\it outer}.

Assign $f$ to a class ${\cal F}_+ \subset {\cal F}$ if it belongs to
$C^\infty \left([0, \infty); {\cal B}\right)$, takes the values in
$\Gamma_0 {\rm Dom\,} A \subset {\cal B}$, and vanishes near $t=0$,
i.e., satisfies ${\rm supp\,} f \subset (0, \infty)$. Also, note
that $f \in {\cal F}_+$ implies $\Pi \left(f(\,\cdot\,)\right) \in
{\cal M}_{{\cal D}}$ (see item 14).
\begin{lemma}\label{L3}
For $f \in {\cal F}_+$, the classical solution $u^f$ to problem
$(\ref{0.1})-(\ref{0.3})$ is represented in the form
\begin{equation}\label{u^f=h-int}
u^f(t)\,=\,h(t)-\int_0^t L^{-\frac{1}{2}} \sin\left[(t-s)
L^{\frac{1}{2}}\right]\, h^{\prime \prime}(s)\,ds\,, \qquad t
\geqslant 0
\end{equation}
with $h:=\Pi \left(f(\,\cdot\,)\right) \in {\cal M}_{{\cal D}}$.
\end{lemma}
{\bf Proof.}\,\,\,Introducing a new unknown
$w=w^f(t):=u^f(t)-\Pi\left(f(t)\right)$ and taking into account C1
(item 17), we easily get the system
\begin{align*}
& w_{tt}+A w=-\Pi\left(f_{tt}(t)\right) &&{\rm in}\,\, {\cal H}, \quad 0<t<\infty   \\
& w|_{t=0}=w_t|_{t=0}=0 &&{\rm in}\,\, {\cal H} \\
& \Gamma_0 w = 0 &&{\rm in}\,\,{\cal B}, \quad 0\leqslant t
<\infty\,.
\end{align*}
With regard to the definition of the operator $L$ (see the axiom
R2), this problem can be rewritten in the form
\begin{align*}
& w_{tt}+L w\,=\,- h_{tt} &&{\rm in}\,\, {\cal H}, \quad 0<t<\infty   \\
& w|_{t=0}=w_t|_{t=0}=0 &&{\rm in}\,\, {\cal H}
\end{align*}
and then solved by the Duhamel formula
\begin{equation*}
w^f(t)\,=-\int_0^t L^{-\frac{1}{2}} \sin\left[(t-s)
L^{\frac{1}{2}}\right]\, h^{\prime \prime}(s)\,ds\,.
\end{equation*}
Returning back to $u^f=w^f + \Pi f$, we arrive at
(\ref{u^f=h-int}). \,\,$\square$
\bigskip

{\bf 20. Reachable sets.}\,\,The sets
\begin{align}
\notag & {\cal U}^t_+\,:=\,\{u^f(t)\,|\,\, f \in {\cal F}_+\}=
\langle\,{\rm see\,} (\ref{u^f=h-int})\, \rangle\,=\\
\label{reachable_from_boundary} &\biggl\{h(t)-\int_0^t
L^{-\frac{1}{2}} \sin\left[(t-s) L^{\frac{1}{2}}\right]\,
h^{\prime \prime}(s)\,ds\,\biggl |\,\,\,\,h=\Pi
f(\,\cdot\,),\,\,\,f \in {\cal F}_+ \biggr\}\,, \quad t \geqslant
0
\end{align}
are said to be {\it reachable from boundary}.

The Green system, which governs the DSBC, determines the certain
pair $L, {\cal D}$, which in turn determines the family $\{{\cal
U}^t_{\cal D}\}$ by (\ref{reachable_sets_U}). Comparing
(\ref{reachable_sets_U}) with (\ref{reachable_from_boundary}), we
easily conclude that the embedding ${\cal U}^t_+ \subset {\cal
U}^t_{\cal D}$ holds. Moreover, the density properties R3 (item
17) enable one to derive $\overline{{\cal U}^t_+} =
\overline{{\cal U}^t_{\cal D}}, \,\,t \geqslant 0$. It is the
latter relation, which inspires the definition
(\ref{reachable_sets_U}) and motivates the terms `reachable sets',
`boundary nest', etc in the general case (item 14), where neither
boundary value space nor boundary operators are defined.
\bigskip

{\bf 21. Illustration.}\,\,Return to the item 19. The DSBC
(\ref{0.1})--(\ref{0.3}) associated with the Riemannian manifold
is governed by the wave equation and is of the form
\begin{align}
\label{0.1Rim} & u_{tt}-\Delta u=0 &&{\rm in}\,\,\,\Omega \times (0, \infty) \\
\label{0.2Rim} & u|_{t=0}=u_t|_{t=0}=0 &&{\rm in}\,\,\,\Omega \\
\label{0.3Rim} & u|_\Gamma =f(t) &&{\rm for}\,\,\,0\leqslant t <
\infty
\end{align}
with a boundary control $f \in {\cal F} = L_2^{\rm loc}\left((0,
\infty); L_2(\Gamma) \right)$. The solution $u=u^f(x,t)$ describes
a wave, which is initiated by boundary sources and propagates from
the boundary into the manifold with the speed 1. For $f \in {\cal
F}_+=C^\infty \left([0, \infty); C^\infty(\Gamma) \right)$
provided ${\rm supp\,} f \subset (0, \infty)$, the solution $u^f$
is classical.

By the finiteness of the wave propagation speed, at a moment $t$
the waves fill a near-boundary subdomain $\Gamma^t:=\{x \in
\Omega\,|\,\,{\rm dist\,} (x, \Gamma) <t\}$. Correspondingly, the
reachable sets ${\cal U}^t_+$ increase as $t$ grows and the
relation ${\cal U}^t_+ \,\subset {\cal H}\Gamma^t, \,\, t
\geqslant 0 $ holds \footnote{Geometric subspaces ${\cal H} A$ are
defined in item 8.}. Closing in ${\cal H}$, we get $\overline
{{\cal U}^t_{\cal D}} \,\subseteq {\cal H}\Gamma^t, \,\, t
\geqslant 0 $
\bigskip

So, if the pair $L, {\cal D}$ (or, equivalently, the operator
$L_0$) appears in the framework of a Green system, then $\{{\cal
U}^t_{\cal D}\}$ introduced by the general definition
(\ref{reachable_sets_U}) can be imagine as the sets of waves
produced by boundary controls. The question arises: What is the
meaning of the corresponding wave spectrum $\Omega_{L, {\cal D}}$
($=\Omega_{L_0}$)? In a sense, it is the question, which this
paper is written for. The answer (section 3) is that, in generic
cases, $\Omega_{L_0}$ is identical to $\Omega$.
\bigskip

{\bf 22. Boundary controllability.}\,\,Return to the abstract DSBC
(\ref{0.1})--(\ref{0.3}) and define its certain property. Begin
with the following observation. Since the class of controls ${\cal
F}_+$ satisfies $\frac{d^2}{dt^2}{\cal F}_+={\cal F}_+$, the
reachable sets (\ref{reachable_from_boundary}) satisfy $A {\cal
U}^t_+ ={\cal U}^t_+$. Indeed, taking $f \in {\cal F}_+$ we have
\begin{equation}\label{relations}
A u^f(t)=\langle\,{\rm see\,} (\ref{0.1})\,\rangle = - u^f_{tt}(t)
= u^{-f^{\prime \prime}}(t) \in {\cal U}^t_+\,.
\end{equation}
By the same relations, $u^f(t)\,=\,A u^g(t)$ holds with
$g=-(\int_0^t)^2 f \in {\cal F}_+$. Hence, the sets ${\cal U}^t_+$
reduce the operator $A$ and its parts $A|_{{\cal U}^t_+}$ are well
defined.
\smallskip

{\bf Definition 13.}\, The DSBC (\ref{0.1})--(\ref{0.3}) is said
to be {\it controllable from boundary} at the time $t=T$ if
$\overline{A|_{{\cal U}^T_+}}=\overline A$ holds, i.e., one
has\footnote{Below the closure is taken in ${\cal H} \times {\cal
H}$; \,${\rm graph\,} A:=\left\{\{y, A y\}\,|\,\,y \in {\rm Dom\,}
A \right\}$.}
\begin{equation}\label{controllability_def_1}
\overline{\left\{\{u^f(T), A u^f(T)\}\,|\,\,f \in {\cal
F}_+\right\}}=\overline{{\rm graph\,} A}=\langle{\rm see\,\,R1}\rangle={\rm
graph\,}L_0^*.
\end{equation}

 Controllability means two things. First, since $A$ is
densely defined in ${\cal H}$, the equality
(\ref{controllability_def_1}) implies $\overline{{\cal U}^t_+}
={\cal H},\,\, t \geqslant T\,,$ i.e., for large times the
reachable sets become rich enough (dense in ${\cal H}$). Second,
the `wave part' $A|_{{\cal U}^T_+}$ of the operator $A$, which
governs the evolution of the system, represents the operator in
substantial.

In applications to problems in bounded domains, such a property
`ever holds' (typically, for large enough times $T$). In
particular, the system (\ref{0.1Rim})--(\ref{0.3Rim}) is
controllable from boundary for any $T>\underset{x \in \Omega}
\max\, {\rm dist\,}(x, \Gamma)$ \cite{BIP97}, \cite{BIP07}.
\medskip

Let us represent the property (\ref{controllability_def_1}) in the
form appropriate for what follows.

Restrict the system (\ref{0.1})--(\ref{0.3}) on a finite time
interval $[0,T]$. Define the Hilbert space of controls ${\cal F}^T
= L_2\left([0, T]; {\cal B}\right)$ and the corresponding smooth class
${\cal F}^T_+ \subset{\cal F}^T$.

Introduce a {\it control operator} $W^T: {\cal F}^T \to {\cal H},
\,\,{\rm Dom} W^T = {\cal F}^T_+,\, W^T f\,:=\, u^f(T)$. Let
$W^T\,=\,U^T\,|W^T|$ be its {\it polar decomposition}, where
$|W^T|:=\left((W^T)^*W^T\right)^{1 \over 2}$ acts in ${\cal F}^T$, and
$U^T$ is an isometry from $\overline{{\rm Ran\,} |W^T|} \subset {\cal
F}^T$ onto $\overline{{\rm Ran\,} W^T} \subseteq {\cal H}$ (see, e.g.,
\cite{BSol}).
\begin{lemma}\label{L4}
If the DSBC $(\ref{0.1})-(\ref{0.3})$ is controllable at $t=T$
then the relation $\overline{\left\{\{|W^T| f,\, |W^T| (-f^{\prime
\prime})\}\,|\,\,f \in {\cal F}_+\right\}}=(U^T)^*L_0^* U^T$
holds.
\end{lemma}
{\bf Proof.}\,\,Represent (\ref{controllability_def_1}) in the
equivalent form $\overline{\{\{W^T f, W^T (-f^{\prime \prime})\}\,|\,\,f
\in {\cal F}_+\}}\newline={\rm graph\,} L_0^*$. Since $\overline{{\rm Ran}
U^T}=\overline{{\cal U}^T}={\cal H}$, the isometry $U^T$ is a unitary operator.
Applying it to the latter representation, one gets the assertion
of the lemma.\,\,$\square$

As a consequence, we conclude the following.
\begin{proposition}\label{P11}
If the DSBC $(\ref{0.1})-(\ref{0.3})$ is controllable at $t=T$
then the operator $|W^T|$ determines the operator $L_0^*$ up to
unitary equivalence.
\end{proposition}
\bigskip

\noindent{\bf 23. Response operator.}\,\,In the DSBC
(\ref{0.1})--(\ref{0.3}) restricted on $[0,T]$, an `input--output'
correspondence is described by the {\it response operator} $R^T:
{\cal F}^T \to {\cal F}^T, \,\,\,{\rm Dom\,} R = {\cal F}^T_+,
\,\,\left(R^Tf\right)(t)\,:=\,\Gamma_1 \left(u^f(t)\right)\,, \,\,
0\leqslant t \leqslant T$.

{\bf Illustration.}\,The response operator of the DSBC
(\ref{0.1Rim})--(\ref{0.3Rim}) is $R^{T}: f\,\mapsto \partial_\nu
u^f|_{\Gamma \times [0,T]}$.
\smallskip

The key fact of the BC-method is that the operator $R^{2T}$
determines the operator $C^T:=(W^T)^*W^T$ through an explicit
formula \cite{BIP97},\cite{DSBC},\cite{BIP07}.
\begin{proposition}\label{P12}
The representation $C^T=\frac{1}{2}(S^T)^* R^{2T}J^{2T} S^T$
holds, where the operator $S^T: {\cal F}^T \to {\cal F}^{2T}$ extends
controls from $[0,T]$ to $[0,2T]$ by oddness w.r.t. $t=T$,
$J^{2T}: {\cal F}^{2T} \to
{\cal F}^{2T},\,\,(J^{2T}f)(t)=\int_0^tf(s)\,ds$.\end{proposition}
Hence, $R^{2T}$ determines the modulus $|W^T|=(C^T)^\frac{1}{2}$.
By Proposition \ref{P11}, we conclude that $R^{2T}$ determines the
operator $L_0^*$ up to unitary equivalence. Since $L_0=L_0^{**}$,
we arrive at the following basic fact.
\begin{proposition}\label{P13}
If the DSBC $(\ref{0.1})-(\ref{0.3})$ is controllable from
boundary at $t=T$ then its response operator $R^{2T}$ determines
the operator $L_0$ up to unitary equivalence.
\end{proposition}
\bigskip

\noindent{\bf 24. Illustration.}\,\,The system
(\ref{0.1Rim})--(\ref{0.3Rim}) is also controllable from boundary.
Such a property is a partial case of the following general fact.

Return to the system (\ref{dynsystem_1})--(\ref{dynsystem_1}). In
our case, the operator $L$ governing its evolution is the
Dirichlet Laplacian $-\Delta$ (item 18). Fix a set $A \in {{\cal O}^{\rm reg}}$.
The reachable sets ${\cal V}^t_{{\cal H} A}$ consist of the waves produced
by sources supported in $A \subset \Omega$. Since the waves
propagate with unit velocity, the embedding ${\cal V}^t_{{\cal H} A}
\subseteq {\cal H} A^t$ holds evidently. The character of this
embedding is a subject of control theory of hyperbolic PDE.

The principal result is that the relation $\overline {{\cal
V}^t_{{\cal H} A}}= {\cal H} A^t$ is valid for any $A \in {{\cal
O}^{\rm reg}}$ and $t\geqslant 0$. It is derived from the
fundamental Holmgren-John-Tataru uniqueness theorem (see,
e.g.,\cite{BIP97}, \cite{BIP07}). In control theory this property
is referred to as a {\it local controllability} of manifolds. In
notation of item 13, it takes the form: $(I_L {\cal H} A)(t) =
{\cal H} A^t$ holds for any $A \in {{\cal O}^{\rm reg}},
\,t\geqslant 0$. Since ${\cal H} A^t=(iM{\cal H} A)(t)$ by the
definition of metric inflation on ${\mathfrak R}^{\cal H}$ (item
8), we arrive at the following formulation of the local
controllability.
\begin{proposition}\label{P14}
The inflations $I_L$ and $iM$ coincide on the lattice ${{\mathfrak R}^{\cal H}}$.
\end{proposition}
\smallskip

Return to the system (\ref{0.1Rim})--(\ref{0.3Rim}) and the
embedding $\overline {{\cal U}^t_{\cal D}}\subseteq {\cal H}
\Gamma^t$ (item 21). The same HJT-theorem implies the equality
$\overline {{\cal U}^t_{\cal D}}= {\cal H} \Gamma^t,\, t\geqslant
0$, which is referred to as a {\it local boundary controllability}
of the manifold $\Omega$.

Recall that the boundary nest ${\mathfrak
u}_{L_0}=\{\overline{{\cal U}^t_{\cal D}}\}_{t \geqslant 0}$
(${\cal D}={\rm Ker} L_0^*$) is defined in item 15. Let
${\mathfrak b}=\{\Gamma^t\}_{t \geqslant 0} \subset {{\cal O}^{\rm
reg}}$ be the family of metric neighborhoods of the boundary
$\Gamma$. Denote $[{\mathfrak b}]=\{[\Gamma^t]\}_{t \geqslant 0}
\subset {\mathfrak R}$ (items 7,8). Boundary controllability of
$\Omega$ is equivalent to the following.
\begin{proposition}\label{P15}
The relation $i[\Gamma^T]=\overline{{\cal U}^t_{\cal D}},\, t
\geqslant 0$ holds. Hence, $i[{\mathfrak b}]={\mathfrak u}_{L_0}$.
\end{proposition}
\smallskip

Boundary controllability implies the following. Since the family
$\{\Gamma^t\}$ exhausts $\Omega$ for any $T \geqslant T_*:=\sup_{x
\in \Omega} {\rm d}(x, \Gamma)$, the boundary nest
$\{\overline{{\cal U}^t_{\cal D}}\}_{t \leqslant T}$ exhausts the
space ${\cal H}$ as $T \geqslant T_*$. By this, the system
(\ref{0.1Rim})--(\ref{0.3Rim}) turns out to be controllable as $T
\geqslant T_*$ \cite{BIP97},\cite{BIP07}.
\smallskip

Hence, by Proposition \ref{P13}, given for a fixed $T\geqslant
2T_*$ the response operator $R^T$ of the system
(\ref{0.1Rim})--(\ref{0.3Rim}) determines the minimal Laplacian
$L_0$ up to unitary equivalence.

\subsection{Stationary DSBC}
{\bf 25. Weyl function.}\,\,Here we follow the paper \cite{Ryzh},
and deal with the same Green system $\{ {\cal H},{\cal B};
A,\Gamma_0, \Gamma_1\}$ and the associated operators $L_0, L$
(item 17).
\smallskip

The problem
\begin{align}
\label{0.1Stat} & \left(A - z \mathbb I\right) u=0 &&{\rm in}\,\,
{\cal H},\,\,\, z \in \mathbb C  \\
\label{0.2Stat} & \Gamma_0 u=\varphi  &&{\rm in}\,\,{\cal B}
\end{align}
is referred to as a {\it stationary} DSBC.  For $\varphi \in
\Gamma_0 {\rm Dom\,} A$ and $z \in {\mathbb C\,}\backslash {\,\rm
spec\,} L$, such a problem has a unique solution $u=u^\varphi_z$,
which is a ${\rm Dom\,} A\,$-valued function of $z$.

The `input--output' correspondence in the system
(\ref{0.1Stat})--(\ref{0.2Stat}) is realized by an operator-valued
function $W(z): {\cal B} \to {\cal B}, \,W(z)\varphi:=\Gamma_1
u^\varphi_z\,\, (z \notin {\,\rm spec\,} L)$. It is called the
{\it Weyl function} and plays the role of data in frequency domain
inverse problems.

The following important fact is established in \cite{Ryzh}. Recall
that a symmetric operator in ${\cal H}$ is said to be {\it
completely non-selfadjoint} if there is no subspace in ${\cal H}$,
in which the operator induces a self-adjoint part.
\begin{proposition}\label{P16}
If the Green system is such that the operator $L_0$ is completely
non-selfadjoint, then the Weyl function determines the operator
$L_0$ up to unitary equivalence.
\end{proposition}
\bigskip

{\bf 26. Illustration.}\,\,Return to item 17. The DSBC
(\ref{0.1Stat})--(\ref{0.2Stat}) associated with the Riemannian
manifold is
\begin{align}
\label{0.1StatRim} & \left(A + z \right) u=0 \qquad {\rm in}\,\,\Omega\\
\label{0.2StatRim} & u|_\Gamma\,=\,\varphi\,,
\end{align}
where $A=-\Delta|_{H^2(\Omega)}$.
\begin{lemma}
The operator $L_0=-\Delta|_{H^2_0(\Omega)}$ is completely
non-selfadjoint.
\end{lemma}
{\bf Proof.}\,\,Assume that there exists a subspace ${\cal K}
\subset {\cal H}$ such that the operator $L_0^{{\cal
K}}:=-\Delta|_{{\cal K} \cap H^2_0(\Omega)}\not= \mathbb O$ is
self-adjoint in ${\cal K}$. In the mean time, $L_0^{{\cal K}}$ is
a part of $L$, which is a self-adjoint operator with the discrete
spectrum. Hence, ${\rm spec}\,L^{{\cal K}}_0$ is also discrete;
each of its eigenfunctions satisfies $-\Delta \phi=\lambda \phi$
in $\Omega$ and belongs to $H^2_0(\Omega)$. The latter implies
$\phi =
\partial_\nu \phi=0$ on $\Gamma$. This leads to $\phi
\equiv 0$ by the well-known E.Landis uniqueness theorem for
solutions to the Cauchy problem for elliptic equations. Hence,
$L^{{\cal K}}_0=\mathbb O$ in contradiction to the
assumption.\,\,$\square$
\smallskip

The Weyl function of the system  is $W(z) \varphi\,=\,
\partial_\nu u^\varphi_z |_\Gamma \,\,(z
\not \in {\rm spec\,} L)$. By the aforesaid, the function $W$
determines the minimal Laplacian $L_0$ of the manifold $\Omega$
up to unitary equivalence.
\bigskip

{\bf 27. Spectral data.}\,\, Besides the Weyl function, there is
one more kind of boundary inverse boundary data associated with
the DSBC (\ref{0.1StatRim})--(\ref{0.2StatRim}). Let
$\{\lambda_k\}_{k=1}^\infty : \,\,\,0<\lambda_1<\lambda_2\leqslant
\lambda_3 \leqslant \dots \to \infty$ be the spectrum of the
Dirichlet Laplacian $L$. Let $\{\phi_k\}_{k=1}^\infty :
\,\,\,L\phi_k=\lambda_k\phi_k$ be the corresponding eigen basis in
${\cal H}$ normalized by $(\phi_k, \phi_l)=\delta_{kl}$.

The set of pairs $\Sigma_\Omega\,:=\,\left\{\lambda_k
;\,\partial_\nu \phi_k |_\Gamma\right\}_{k=1}^\infty$ is called
the (Dirichlet) {\it spectral data} of the manifold $\Omega$.

The well-known fact is that these data determine the Weyl function
and vice versa (see, e.g., \cite{Ryzh}). Hence, $\Sigma_\Omega$
determines the minimal Laplacian $L_0$ up to unitary equivalence.
However, such a determination can be realized not through $W$ but
in more explicit way.

Namely, let $U: {\cal H} \to \widetilde {\cal H} := \emph{l}_2$,
$U y\,=\,\widetilde y:=\{(y,\phi_k)\}_{k=1}^\infty$ be the Fourier
transform that diagonalizes $L$: $ \widetilde L\,:=\,U L U^*={\rm
diag\,}\{\lambda_1,\,\lambda_2\,,\, \dots\}$. For any harmonic
function $a \in {\cal A}$, its coefficients are
$(a,\phi_k)\,=\,-\,\frac{1}{\lambda_k} \int_\Gamma a\,\partial_\nu
\phi_k\,d\Gamma$ that can be verified by integration by parts.
Therefore, the spectral data $\Sigma_\Omega$ determine the image
$\widetilde {\cal A}:=U {\cal A} \subset \widetilde {\cal H}$ and
its closure $\widetilde {\cal D}=U {\cal D}=\overline{\widetilde {\cal
A}}$. Thus, the determination $\Sigma_\Omega\,\Rightarrow
\widetilde L\,, \widetilde {\cal D}$ occurs.

In the mean time, the relation C3 (item 17) implies $\widetilde
L_0\,=\,U^* L_0 U\,=\, \widetilde L |_{{\widetilde
L}^{-1}\left[\widetilde {\cal H} \ominus \widetilde {\cal
D}\right]}$ by isometry of $U$. Thus, $\widetilde L_0$ is a
unitary copy of $L_0$ constructed via the spectral data.

\section{Reconstruction of manifolds}

\subsection{Inverse problems}
{\bf 28. Setup.}\,\,In inverse problems (IP) for DSBC associated
with manifolds, one needs to recover the manifold via its boundary
inverse data\footnote{In concrete applications (acoustics,
geophysics, electrodynamics, etc), these data formalize the
measurements implemented at the boundary.}. Namely,
\smallskip

\noindent {\bf IP 1.}\,\,\,given for a fixed $T>2\underset{x \in
\Omega} \max\, {\rm dist\,}(x, \Gamma)$ the response operator
$R^{T}$ of the system (\ref{0.1Rim})--(\ref{0.3Rim}), to recover
the manifold $\Omega$
\smallskip

\noindent {\bf IP 2.}\,\,\,given the Weyl function $W$ of the
system (\ref{0.1StatRim})--(\ref{0.2StatRim}), to recover the
manifold $\Omega$
\smallskip

\noindent {\bf IP 3.}\,\,\,given the spectral data
$\Sigma_\Omega$, to recover the manifold $\Omega$.
\smallskip

\noindent The problems are called {\it time-domain}, {\it
frequency-domain}, and {\it spectral} respectively.

Setting the goal to determine an unknown manifold from its
boundary inverse data, we have to keep in mind the evident
nonuiqueness of such a determination: all {\it isometric}
manifolds with the mutual boundary have the same data. Therefore,
the only reasonable understanding of `to recover' is to construct
a manifold, which possesses the prescribed data \cite{BIP07}.

As we saw, the common feature of problems IP 1--3 is that their
data determine the minimal Laplacian $L_0$ up to unitary
equivalence. By this, each kind of data determines the wave
spectrum $\Omega_{L_0}$ up to isometry (see Proposition
\ref{P10}). As will be shown, for a wide class of manifolds the
relation $\Omega_{L_0} \overset{\rm isom}= \Omega$ holds. Hence,
for such manifolds, to solve the IPs it suffices to extract a
unitary copy $\tilde L_0$ from the data, find its wave spectrum
$\Omega_{\tilde L_0} \overset{\rm isom}=\Omega_{L_0}$, and thus
{\it get an isometric copy of} $\Omega$. It is the program for the
rest of the paper.
\bigskip

{\bf 29. Simple manifolds.}\,\, Recall that we deal with a compact
smooth Riemannian manifold $\Omega$ with the boundary $\Gamma$.
The family ${\mathfrak b}=\{\Gamma^t\}_{t \geqslant 0}$ consists
of metric neighborhoods of $\Gamma$. Nets and dense lattices were
introduced in item 9. ${\mathfrak L}[M,{\mathfrak b}] \subset
{{\cal O}^{\rm reg}}$ is the minimal $M$-invariant (sub)lattice,
which contains ${\mathfrak b}$.
\smallskip

We say $\Omega$ to be a {\it simple manifold} if the lattice
${\mathfrak L}[M,{\mathfrak b}]$ is dense in ${{\cal O}^{\rm reg}}$.
\smallskip

The evident obstacle for a manifold to be simple is its
symmetries. For a ball $\Omega=\{x \in {\mathbb
R}^n\,|\,\,|x|\leqslant 1\}$, the lattice ${\mathfrak
L}[{\mathfrak b},M]$ consists of sums of `annuluses' of the form
$\{x \in \Omega\,|\,\,0\leqslant a<|x|<b\leqslant 1\}$. Surely,
such a system is not a net in the ball. A plane triangle is simple
if and only if its legs are pair-wise nonequal. Easily checkable
sufficient conditions on the shape of $\Omega \subset {\mathbb
R}^n$, which provide the simplicity, are proposed in
\cite{BKac89}. They are also appropriate for Riemannian manifolds
and show that simplicity is a generic property: it can be provided
by arbitrarily small smooth variations of the boundary $\Gamma$
\footnote{Presumably, any compact manifold with trivial symmetry
group is simple but it is a conjecture. In the mean time, for
noncompact manifolds this is not true.}.
\bigskip

{\bf 30. Solving IPs.}\,\,The following result provides
reconstruction of $\Omega$.
\begin{theorem}\label{Theorem_1}
Let $\Omega$ be a simple manifold, $L_0=-\Delta|_{H^2_0(\Omega)}$
the minimal Laplacian, $\Omega_{L_0}$ its wave spectrum. There
exists an isometry (of metric spaces) $i_*$ that maps
$\Omega_{L_0}$ onto $\Omega$, the relation $i_*(\partial
\Omega_{L_0})=\Gamma$ being valid.
\end{theorem}
{\bf Proof.}\,\, Denote $[{\mathfrak b}]:= \{[\Gamma^t]\}_{t
\geqslant 0} \subset {\mathfrak R}$. Let ${\mathfrak
L}\left[M,[{\mathfrak b}]\right]\subset {\mathfrak R}$ be the
image of ${\mathfrak L}[M,{\mathfrak b}]$ through the `projection'
$A \mapsto [A]$ (item 7).

Propositions \ref{P14},\ref{P15} imply $i {\mathfrak L}\left[M,[{\mathfrak b}]\right]=
{\mathfrak L}\left[iM, i[{\mathfrak b}]\right]={\mathfrak L}[I_L, {\mathfrak u}_{L_0}]={\mathfrak L}_{L_0}\subset
{{\mathfrak R}^{\cal H}}$.

Taking into account the simplicity condition and applying
Proposition \ref{P9} to the case ${\cal N}={\mathfrak L}[M,{\mathfrak b}]$, we conclude
that $\Omega_{L_0}$ is isometric to $(\Omega, {\rm d})$. The
isometry is realized by the bijection $i_*: i[\alpha] \mapsto
x_\alpha$.
\smallskip

To compare the atoms $i[\alpha]$, which constitute $\Omega_{L_0}$,
with the boundary nest ${\mathfrak u}_{L_0}$ is in fact to compare
the metric neighborhoods $\{x_\alpha\}^t$ with the metric
neighborhoods $\{\Gamma^t\}$. Since $\{x_\alpha\}^t \subset
\Gamma^t, \,\,t \geqslant 0$ is valid if and only if $x_\alpha \in
\Gamma$, we conclude that $i_*(\partial \Omega_{L_0})=
\Gamma$.\,\,\,$\square$
\medskip

Thus, to solve the IPs 1-3 in the case of simple $\Omega$, it
suffices to determine (from the inverse data) a relevant unitary
copy $\tilde L_0$ of the minimal Laplacian, and then find its wave
spectrum $\Omega_{\tilde L_0}$.
\bigskip

{\bf 31. Remarks.}\,\,

{\bf 1.}\,\,Regarding non-simple manifolds, note the following. If
the
 symmetry group of $\Omega$ is nontrivial then, presumably,
$\Omega_{L_0}$ is isometric to the properly metrized set of the
group orbits. Such a conjecture is motivated by the following
easily verifiable examples.
\smallskip

\noindent $\bullet$\,\,For a ball $\Omega=\{x \in {\mathbb
R}^n\,|\,\,|x| \leqslant r\}$, the spectrum $\Omega_{L_0}$ is
isometric to the segment $[0,r] \subset \mathbb R$. Its boundary
$\partial \Omega_{L_0}$ is identical to the endpoint $\{0\}$.
\smallskip

\noindent $\bullet$\,\,For an ellipse $\Omega=\{(x, y)\in {\mathbb
R}^2\,|\,\,\frac{x^2}{a^2} + \frac{y^2}{b^2} \leqslant 1\}$,
$\Omega_{L_0}$ is isometric to its quarter $\Omega \cap \{(x,
y)\,|\,\,x \geqslant 0,\, y \geqslant 0 \}$, whereas $\partial
\Omega_{L_0}$ is isometric to $\{(x, y)\,|\,\,\frac{x^2}{a^2} +
\frac{y^2}{b^2} =1,\,\,x \geqslant 0,\, y \geqslant 0 \}$.
\smallskip

\noindent $\bullet$\,\,Let $\omega \subset \left\{(x_1, x_2)\in
{\mathbb R}^2\,|\,\,x_2>0\right\}$ be a compact domain with the
smooth boundary. Let $\Omega$ be a torus in ${\mathbb R}^3$, which
appears as result of rotation of $\omega$ around the $x_1$-axis.
Then $\Omega_{L_0}\overset{\rm isom}=\omega$ and $\partial
\Omega_{L_0} \overset{\rm isom}=\partial \omega$.
\medskip

{\bf 2.}\,\,In applications a possible lack of simplicity is not
an obstacle for solving problems IP 1--3 because their data not
only determine a copy of $L_0$ but contain substantially more
information about $\Omega$. Roughly speaking, the matter is as
follows. When we deal with these problems, the boundary $\Gamma$
is given. By this, instead of the boundary nest ${\mathfrak
u}_{L_0}$ of the sets reachable from {\it the whole} $\Gamma$ (see
(\ref{reachable_from_boundary})), we can use the much richer
family ${\mathfrak u}^\prime_{L_0}=\{{\cal U}^t_\sigma\}_{t
\geqslant 0, \,\sigma \subset \Gamma}$ of sets reachable from {\it
any patch} $\sigma \subset \Gamma$ of positive measure
\footnote{More precisely, ${\cal U}^t_\sigma$ consists of the
solutions (waves) $u^f(t)$ produced by the boundary controls $f$
supported on $\sigma \times [0, \infty)$}. Therefore, even though
the density of the lattice ${\mathfrak L}[I_L,{\mathfrak u}_{L_0}
]$ in ${{\mathfrak R}^{\cal H}}$ may be violated by symmetries,
the lattice ${\mathfrak L}[I_L,{\mathfrak u}^\prime_{L_0}]$ is
always dense. As a result, the wave spectrum corresponding to the
dense lattice turns out to be isometric to $\Omega$. The latter is
the key fact, which enables one to reconstruct $\Omega$: see
\cite{BSobolev} for detail.
\medskip

{\bf 3.}\,\,The spectra $\Omega^{\rm nest}_{L_0}$ and
$\Omega^{\rm eik}_{L_0}$ are also appropriate for reconstruction. If
$\Omega$ is simple, one has $\Omega_{L_0}\overset{\rm
isom}=\Omega^{\rm nest}_{L_0}\overset{\rm
isom}=\Omega^{\rm eik}_{L_0}\overset{\rm isom}=(\Omega, {\rm d})$
\cite{BSobolev},\cite{BArXiv}.
\medskip

{\bf 4.}\,\,If $\Omega$ is noncompact, the definition of
simplicity remains to be meaningful, local controllability is in
force, and ${\cal H}=\cup_{t>0}\overline{{\cal U}^t_{\cal D}}$ holds. One can show that
the response operator $R^{T}$ known {\it for all} $T>0$ determines
the simple manifold up to isometry. Also, defining mutatis
mutandis the Weyl function and spectral data for a noncompact
$\Omega$, one can obtain the same result: these data determine the
simple manifold up to isometry.
\bigskip

{\bf 32. Algebras in reconstruction.}\,\,Recall that the von
Neumann algebra ${\mathfrak N}_{\mathfrak L} \subset {{\mathfrak B}(\cal H)}$ associated with the lattice
${\mathfrak L} \subset {{\mathfrak L}(\cal H)}$ was introduced in item 11. In the bounded case,
along with ${\mathfrak N}_{\mathfrak L} $ one can define the algebra ${\mathfrak C}_{\mathfrak L}$ as the
minimal norm-closed subalgebra of ${{\mathfrak B}(\cal H)}$, which contains all
maximal eikonals.

For the algebras ${\mathfrak N}_{L_0};={\mathfrak N}_{{\mathfrak L}_{L_0}}$ and
${\mathfrak C}_{L_0};={\mathfrak C}_{{\mathfrak L}_{L_0}}$ associated with manifold, the
following holds \cite{BSobolev},\cite{BArXiv}.
\smallskip

\noindent {(i)}\,\,Both of these algebras are commutative. The
embedding ${\mathfrak C}_{L_0} \subset {\mathfrak N}_{L_0}$ is dense in the strong
operator topology in ${{\mathfrak B}(\cal H)}$.
\smallskip

\noindent {(ii)}\,\,If $\Omega$ is simple then ${\mathfrak C}_{L_0}$ is
isometrically isomorphic to the algebra $C(\Omega)$ of continuous
functions. By this, its spectrum\footnote{i.e., the set of maximal
ideals of ${\mathfrak C}_{L_0}$ \cite{Mur}.} $\hat{\mathfrak C}_{L_0}$ is homeomorphic
to $\Omega$.
\smallskip

These results are applied to reconstruction by the scheme $\{\rm
inverse\, data\} \Rightarrow {\mathfrak C}_{L_0} \Rightarrow \hat{\mathfrak C}_{L_0}
\Rightarrow \Omega$ \cite{BSobolev},\cite{BArXiv}.

Note that commutativity is derived from local controllability of
the system (\ref{0.1Rim})--(\ref{0.3Rim}). In the corresponding
dynamical system on a graph, a lack of controllability occurs and,
as a result, these algebras turn out to be {\it noncommutative}
\footnote{N.Wada, private communication.}. This leads to problems
and difficulties in reconstruction, which are not overcome yet. In
particular, the relations between the spectra $\Omega_{l_0}$ and
$\hat{\mathfrak C}_{L_0}$ are not clear.

\subsection{Comments}
{\bf 33. A look at isospectrality}\,\,\, Let ${\rm
spec}\,L=\{\lambda_k\}_{k=1}^\infty$ be the spectrum of the
Dirichlet Laplacian on $\Omega$ (item 27). The question: "Does
${\rm spec}\,L$ determine $\Omega$ up to isometry?" is a version
of the classical M.Kac's drum problem\cite{Kac}. The negative
answer is well known (see, e.g., \cite{BCDS_isosp_dom}) but, as
far as we know, the satisfactory description of the set of
isospectral manifolds is not obtained yet. The following is some
observations in concern with such a description.

Assume that we deal with a simple $\Omega$. In accordance with
Theorem 1, such a manifold is determined by any unitary copy
$\widetilde L_0$ of the operator $L_0 \subset L$. If the spectrum
of $L$ is given, to get such a copy it suffices to possess the
Fourier image $\widetilde {\cal D} =U{\cal D}$ of the harmonic
subspace in $\widetilde {\cal H} = \emph{l}_2$: see C3, item 17
\footnote{It is the fact, which is exploited in \cite{BKac89}}. In
the mean time, as is evident, if $\Omega$ and $\Omega'$ are
isometric, then the corresponding images are identical:
$\widetilde {\cal D}=\widetilde {\cal D}'$. Therefore, $\Omega$
and $\Omega'$ are isospectral but not isometric if and only if
$\widetilde {\cal D}\not=\widetilde {\cal D}'$. In other words,
the subspace $\widetilde {\cal D}$ is a relevant `index', which
distinguishes the isospectral manifolds.

As an image of harmonic functions, which is admissible for the
given $\widetilde L={\rm diag\,}\{\lambda_1, \lambda_2,
\,\dots\}$, a subspace $\widetilde {\cal D} \subset \emph{l}_2$
has to obey the following conditions:
\smallskip

\noindent {\bf 1.}\,\, A lineal set ${\cal L}_{\widetilde {\cal
D}} :={\widetilde L}^{-1}\left[\emph{l}_2 \ominus \widetilde {\cal
D}\right]$ is dense in $\emph{l}_2$, whereas replacement of
$\widetilde {\cal D}$ by any wider subspace ${\widetilde{\cal
D}}^\prime \supset \widetilde {\cal D}$ leads to the lack of
density: ${\rm clos\,} {\cal L}_{\widetilde {\cal
D}^\prime}\not=\emph{l}_{2}$
\smallskip

\noindent {\bf 2.}\,\,Extending an operator $\widetilde L|_{{\cal
L}_{\widetilde {\cal D}}}$ by Friedrichs, one gets $\widetilde L$.

In the mean time, taking {\it any} subspace $\widetilde {\cal D}
\subset \emph{l}_2$ obeying 1,2 \footnote{such subspaces do exist
(M.M.Faddeev, private communication)}, one can construct a
symmetric operator $\widetilde L_0$ by C3, and then find its wave
spectrum $\Omega_{\widetilde L_0}$ as a candidate to be a drum.
However, the open question is whether such a `drum' is human (is a
manifold).
\bigskip

\noindent{\bf 34. Wave model.}\,\,Return to the abstract system
(\ref{0.1})--(\ref{0.3}) and assume it to be controllable at
$t=T$. Reduce the system to the interval $0\leqslant t \leqslant
T$. Recall that the image and control operators ${\cal I}: {\cal
H} \to {\cal G}$ and $W^T: {\cal F}^T \to {\cal H}$ were
introduced in items 16 and 22 respectively. The composition
$V^T:={\cal I} W^T: {\cal F}^T \to {\cal G}$ is called a {\it
visualizing operator} \cite{BIP97}, \cite{DSBC}, \cite{BIP07}.

Let the response operator $R^{2T}$ be given. The following is a
way to construct a canonical `functional' model of the operator
$L_0^*$.
\smallskip

\noindent $\bullet$\,\,$R^{2T}$ determines the operator $|W^T|$ in
${\cal F}^T$  (item 23). In what follows, it is regarded as a {\it
model control operator} $\tilde W^T :=|W^T|$, which acts from
${\cal F}^T$ to a {\it model inner space} $\tilde {\cal H} :={\cal F}^T$.
\smallskip

\noindent $\bullet$\,\,Determine the operator ${\tilde L_0}^*$ in
$\tilde {\cal H}$ as the operator of the graph
\newline$\overline{\left\{\{\tilde W^T f,\, \tilde W^T(-f^{\prime
\prime})\}\,|\,\,f \in {\cal F}_+\right\}}$ (Lemma \ref{L4}, item
22). Find $\tilde L_0 ={\tilde L_0}^{**}$.
\smallskip

\noindent $\bullet$\,\,Find the wave spectrum $\Omega_{\tilde
L_0}$ and recover the germ space $\tilde {\cal G}$ on it. Determine the
image operator $\tilde {\cal I}: \tilde {\cal H} \to \tilde {\cal G}$. Compose
the visualizing operator $\tilde V^T =\tilde {\cal I}\,\tilde W^T:
{\cal F}^T \to \tilde{\cal G}$.
\smallskip

\noindent $\bullet$\,\,Define $({L_0^{\rm mod}})^*$ as an operator in
$\tilde {\cal G}$ determined by the graph \newline $\left\{\{\tilde V^T
f,\, \tilde V^T (-f^{\prime \prime})\}\,|\,\,f \in {\cal
F}_+\right\}$.
\smallskip

Surely, it is just a draft\footnote{Some detail see in
\cite{BArXiv}, sec 3.4 .} of the model and plan for future work:
one needs to endow the germ space ${\cal G}$ with relevant Hilbert
space attributes. Presumably, in `good cases', ${\cal G}=L_{2,
\mu}(\Omega_{L_0})$. Also, the model operator is expected to be
{\it local}: ${\rm supp} (L_0^{\rm mod})^* y \subseteq {\rm supp}
y$, whereas the model trace operators ${\tilde \Gamma}_{0,1}$ are
connected with the restriction $y \mapsto y|_{\partial
\Omega_{L_0}}$\footnote{As far as we know, the known models of
symmetric operators do not possess such properties
\cite{Shtraus}.}. Hopefully, the collection $\{\tilde{\cal
G},{\cal B}; (L_0^{\rm mod})^* ,{\tilde \Gamma}_0, {\tilde
\Gamma}_1\}$ constitutes the Green system, which is a canonical
model of the original $\{ {\cal H},{\cal B}; A,\Gamma_0,
\Gamma_1\}$. The model is determined by $R^{2T}$.

Such a model is in the spirit of general system theory \cite{KFA},
where it would be regarded as a {\it realization} relevant to the
{\it transfer operator function} $R^{2T}$. Remarkable point is the
role of a time in its construction.
\bigskip

\noindent{\bf 35. Open question.}\,\,\,For any operator $L_0$ of
the class under consideration, the lattice ${\mathfrak L}_{L_0}$
is a well-defined object, ${\mathfrak L}_{L_0}\not=\{0\}$ being
hold. We have neither a proof nor a counterexample to the
following principal conjecture: $\Omega_{L_0}\not=\emptyset$.
However, there is example of the operator $L_0$ such that
$\Omega_{L_0}$ consists of a single point.
\bigskip

\noindent{\bf 36. A bit of philosophy.}\,\,\,In applications, the
external observer pursues the goal to recover a manifold $\Omega$
via measurements at its boundary $\Gamma$. The observer prospects
$\Omega$ with waves $u^f$ produced by boundary controls. These
waves propagate into the manifold, interact with its inner
structure and accumulate information about the latter. The result
of interaction is also recorded at $\Gamma$. The observer has to
extract the information about $\Omega$ from the recorded.

By the rule of game in IPs, the manifold itself is invisible
(unreachable) in principle. Therefore, the only thing the observer
can hope for, is to construct from the measurements an {\it image}
of $\Omega$ possibly resembling the original. By the same rule,
the only admissible material for constructing is the waves $u^f$.
To be properly formalized, such a look at the problem needs two
things:
\smallskip

\noindent $\bullet$ an object that codes exhausting information
about $\Omega$ and, in the mean time, is determined by the
measurements
\smallskip

\noindent $\bullet$ a mechanism that decodes this information.

\noindent Resuming our paper, the first is the minimal Laplacian
$L_0$, whereas to decode information is to determine its wave
spectrum constructed from the waves $u^f$. It is $\Omega_{L_0}$,
which is a relevant image of $\Omega$.
\medskip

The given paper promotes an algebraic trend in the BC-method
\cite{BSobolev}, by which {\it to solve IPs is to find spectra of
relevant lattices and algebras}. An attempt to apply this
philosophy to solving new problems would be quite reasonable. An
encouraging fact is that in all above-mentioned unsolved IPs of
anisotropic elasticity and electrodynamics, graphs, etc, the wave
spectrum $\Omega_{L_0}$ does exist. However, to recognize how it
looks like and verify (if true!) that $\Omega_{L_0}$ is isometric
(homeomorphic) to $\Omega$ is difficult in view of very
complicated structure of the corresponding reachable sets ${\cal
U}^t$.

\subsection{Appendix}
{\bf 37. Basic lemma.}\,\,Recall the notation: for a set $A
\subset \Omega$, $\overline A$ is its metric closure, ${\rm int} A$ is the
set of interior points, $A^t$ is the metric neighborhood of radius
$t$, $A^0:=A$. If $A \in {\cal O}$ then $A\subseteq {\rm int} \overline A$ and
$\overline A=\overline{{\rm int} \overline A}$ holds.

Return to item 5. Let $f=f(t), \, t\geqslant 0$ be an element of
$\overline{M{\cal O}}$. Define the set $\dot f:= \cap_{t>0}f(t)
\subset \Omega$. Define the functions $f_*(t)=(M\dot f)(t)={\dot
f}^t$ as $t>0,\,\,f_*(0)=f(0)$ and $f^*(t)={\rm int}
\overline{{\dot f}^t},\,\, t\geqslant 0$.
\begin{lemma}\label{Appendix}
(i)\,\,If $f \not=0_{{\cal F}_{\cal O}}$ then $\dot
f=\overline{\dot f}\not= \emptyset$ and the relations $f_*
\leqslant f \leqslant f^*$ hold in ${\cal F}_{\cal O}$. (ii)\,\,If
$f$ and $g$ satisfy $\dot f= \dot g$ then
$\overline{f(t)}=\overline{g(t)}\,\,(=\overline{{\dot f}^t})$ as
$t \geqslant 0$.
\end{lemma}
{\bf Proof.}\,\,

\noindent{\bf 1.}\,\,If $f\nwarrow f_j \in M{\cal O}$ then
$f(t)=\cup_{j\geqslant 1}f_j(t),\,t \geqslant 0$. Therefore,
$f_k(0)\subseteq \cup_{j\geqslant 1}f_j(0) \subseteq
\cup_{j\geqslant 1}(f_j(0))^t=\cup_{j\geqslant 1}f_j(t)=f(t), \,t
\geqslant 0$. Hence, $\dot f \supseteq f_k(0)\not= \emptyset$.
\medskip

\noindent{\bf 2.}\,\,If $f\swarrow f_j \in M{\cal O}$ then
$f(t)={\rm int}\cap_{j \geqslant 1}f_j(t), \,t \geqslant 0$.
Define a closed set $F=\cap_{j\geqslant 1}\overline{f_j(0)}
\subset \Omega$ and show that $F\not= \emptyset$.

Assume $F=\emptyset$. Since $\overline{f_{j+1}(0)}\subseteq
\overline{f_{j}(0)}$, for any $x\in \Omega$ and $t>0$ there is
$j_0=j_0(x,t)$ such that $\overline{\{x\}^t} \cap
\overline{f_j(0)}$ as $j > j_0$. Indeed, otherwise, by assumptions
A1,2, the ball $\overline{\{x\}^t}$ has to contain the points of
$F$. Hence, $x \notin (\overline{f_j(0)})^t$ as $j>j_0$. Since $x$
is arbitrary, we have $\emptyset=\cup_{j\geqslant 1}(\overline
{f_j(0)})^t=\cup_{j\geqslant 1}(f_j(0))^t=\cup_{j\geqslant
1}f_j(t)$. Therefore $f(t)={\rm int} \cap_{j\geqslant
1}f_j(t)=\emptyset$, i.e., $f(t)=0_{\cal O},\, t\geqslant 0$. It
means that $f=0_{{\cal F}_{\cal O}}$ in contradiction with
assumptions of the lemma. So, $F\not=\emptyset$.
\medskip

\noindent{\bf 3.}\,\,Show that $F=\dot f$, i.e., $F$ does not
depend on $\{f_j\}$\footnote{However, the limit $f$ can depend on
$\{f_j\}$: there are examples for $\Omega={\mathbb R}^n$!}. For
every $j \geqslant 1$, we have
$\overline{f_j(0)}=\cap_{t>0}(f_j(0))^t=\cap_{t>0}f_j(t) \supseteq
\cap_{t>0}f(t)=\dot f$. Hence $F=\cap_{j\geqslant
1}\overline{f_j(0)} \supseteq \dot f$.

On the other hand, the monotonicity $\overline{f_{j+1}(0)}
\subseteq \overline{f_{j}(0)}$ implies \newline
$F=\cap_{j\geqslant 1}\overline{f_j(0)} \subseteq
(\cap_{j\geqslant 1}\overline{f_j(0)})^t\subseteq \cap_{j\geqslant
1}(\overline{f_j(0)})^t$. Since the next to the last set is open
as $t>0$, we have $F \subseteq (\cap_{j\geqslant
1}\overline{f_j(0)})^t\subseteq {\rm int} \cap_{j\geqslant
1}(\overline{f_j(0)})^t \subseteq {\rm int} \cap_{j\geqslant
1}({f_j(0)})^t\subseteq {\rm int} \cap_{j\geqslant 1}f_j(t)=f(t)$
for all $t>0$. Hence $F \subseteq \cap_{t>0}f_j(t)=\dot f$, and we
arrive at $F=\dot f$.

Thus, we obtain $F=\dot f \not=\emptyset$.
\medskip

\noindent{\bf 4.}\,\,Show that $f_* \leqslant f$. Choosing $M{\cal
O} \ni f_j\searrow f$, for $t>0$ one has ${\dot f}^t=F^t \subseteq
(\overline{f_j(0)})^t=({f_j(0)})^t=f_j(t)$. This implies ${\dot
f}^t \subseteq \cap_{j\geqslant 1}f_j(t)$. Since ${\dot f}^t$ is
an open set, the embedding ${\dot f}^t \subseteq {\rm int}
\cup_{j\geqslant 1}f_j(t)=f(t)$ holds. The latter means that
$f_*(t) \leqslant f(t)$ in ${\cal O}$ as $t>0$. The definition of
$f_*$ at $t=0$ leads to $f_*(t) \leqslant f(t),\,t>0$ in ${\cal
O}$, i.e., $f_* \leqslant f$ in ${\cal F}_{\cal O}$.

Show that $f\leqslant f^*$. Choose $M{\cal O} \ni f_j\searrow f$
that means $f(t)={\rm int}\cap_{j \geqslant 1}f_j(t), \,t
\geqslant 0$. For $t=0$ one has $f(0)={\rm int}\cap_{j \geqslant
1}f_j(0)\subset {\rm int}\cap_{j \geqslant
1}\overline{f_j(0)}={\rm int} \dot f ={\rm int} \overline{\dot
f}=f^*(0)$. For $t>0$, with regard to monotonicity of
$\{f_j\}\!\!\downarrow$, we have $f(t)={\rm int}\cap_{j \geqslant
1}f_j(t)={\rm int}\cap_{j \geqslant 1}(f_j(0))^t={\rm int}\cap_{j
\geqslant 1}{\overline{(f_j(0))}}^t \subseteq {\rm int}
\overline{(\cap_{j\geqslant 1}\overline{f_j(0)})^t}={\rm int}
\overline{{\dot f}^t}=f^*(t)$. Hence $f\leqslant f^*$ is valid.

Thus, the part (i) of the lemma is proven.
\medskip

\noindent{\bf 5.}\,\,For $t>0$, since ${\dot f}^t$ is an open set,
one has $\overline{{\dot f}^t}=\overline{{\rm int}{\overline{{\dot f}^t}}}$. Therefore,
$\overline{f_*(t)}=\overline{f^*(t)}=\overline{{\dot f}^t}$, and (i) implies
$\overline{f(t)}=\overline{{\dot f}^t}$. Hence, $\overline{f(t)}=\overline{{\dot f}^t}
=\overline{{\dot g}^t}=\overline{g(t)}$ as $t>0$.

Let $t=0$. Choosing $M{\cal O} \ni f_j\searrow f$, one has
$f(0)={\rm int} \cap_{j\geqslant 1}f_j(0) \subseteq {\rm int}
\cap_{j\geqslant 1}\overline{f_j(0)}={\rm int} \dot f$. Hence
$\overline{f(0)} \subseteq \overline{{\rm int} \dot f}\subseteq
\dot f$. Show that $\overline{f(0)} = \overline{{\rm int} \dot
f}$. Indeed, assuming the opposite, one can find $x \in \dot f$
separated from  $\overline{f(0)}$ with a positive distance. In the
mean time, defining $f^\varepsilon$ by
$f^\varepsilon(t)=(\overline {f(0)})^{\varepsilon +t},\ t\geqslant
0$, we get $f_*(t)={\dot f}^t \subseteq f(t) \subset
f^\varepsilon(t)$. However, the relation ${\dot f}^t \subset
f^\varepsilon(t)$ is impossible for small enough $t$ and
$\varepsilon$ by the choice of $x$. Hence, $\overline{f(0)} =
\overline{{\rm int} \dot f}$ does hold.

The latter implies $\overline{f(0)} = \overline{{\rm int} \dot f}=
\overline{{\rm int} \dot g}= \overline{g(0)}$. Thus, we get
$\overline{f(t)}=\overline{g(t)}$ for all $t \geqslant 0$ and
prove (ii). \,\,\,$\square$

\bigskip

\noindent{\bf Key words:} symmetric semi-bounded operator, lattice
with inflation, evolutionary dynamical system, wave spectrum,
reconstruction of manifolds
\smallskip

\noindent{\bf MSC:} 47XX, 47A46, 46Jxx, 35R30, 06B35.

\end{document}